\documentclass[10pt,final,fleqn]{article}
\usepackage{a4}
\usepackage{amsmath}%
\usepackage{amstext}%
\usepackage{amssymb}%
\usepackage{amsthm}
\usepackage{comment}
\usepackage{listings}

\usepackage[final]{graphicx}
\usepackage{subcaption}

\usepackage[margin=0.5in]{geometry}

\usepackage{graphicx}
\usepackage{epstopdf}

\usepackage{xcolor}

\newtheorem{theorem}{Theorem}

\newtheorem{axiom}{Axiom}

\newtheorem{corollary}[theorem]{Corollary}

\newtheorem{definition}[axiom]{Definition}

\newtheorem{lemma}[theorem]{Lemma}

\newtheorem{proposition}[theorem]{Proposition}
\newtheorem{rem}[theorem]{Remark}
\newenvironment{remark}{\rem\rm}{\endrem}

\newcommand{\R}{\mathbb{R}}%
\newcommand{\e}{\varepsilon}%
\newcommand{\ol}{\overline}%

\newcommand{\n}{{\nabla}}

\newcommand{\To}{\longrightarrow}
\def\a{\alpha}

\def\b{\beta}
\def\e{\epsilon}

\def\t{\theta}
\def\g{\gamma}

\def\<{\langle}
\def\>{\rangle}

\DeclareMathOperator*\crit{crit}
\DeclareMathOperator*\dist{dist}

\title{Second-order dynamical systems with superlinear convergence rates}

\author{Mikhail Karapetyants \thanks{Institute of Science and Technology Austria (ISTA), Am Campus 1, 3400 Klosterneuburg, Austria} \and Szil\'ard Csaba L\'aszl\'o\thanks{Department of Mathematics, Technical University of Cluj-Napoca, Memorandumului 28, Cluj-Napoca, Romania. E-mail: szilard.laszlo@math.utcluj.ro.}}

\medskip

\begin{document}

\maketitle

\begin{abstract}
This paper investigates a second-order dynamical system with variable coefficients designed for the minimization of smooth non-convex functions. We establish the convergence of the generated trajectories to a critical point of the objective function, provided that a regularization of the objective function satisfies the Kurdyka--\L{}ojasiewicz (KL) property. Furthermore, we derive convergence rates for both the trajectories and the objective function values in terms of the \L{}ojasiewicz exponent. We demonstrate that employing an unbounded damping term within the dynamical system significantly accelerates convergence compared to existing results in the literature. Specifically, we obtain superlinear convergence rates both for the trajectories $x(t)$ toward a critical point and for the potential energy $g(x(t))-\min g$, as $t \to +\infty$. Notably, these results remain valid and directly applicable in the case of strongly convex objective functions.
\end{abstract}



\section{Introduction}

In this paper, we investigate the asymptotic behavior of trajectories generated by the second-order dynamical system with time-dependent coefficients
\begin{equation}\label{dysy}\tag{DS}
\left\{
\begin{aligned}
&\ddot{x}(t) + \alpha(t)\dot{x}(t) + \beta(t)\nabla g(x(t)) = 0,\\
&x(t_0) = u_0, \quad \dot{x}(t_0) = v_0,
\end{aligned}
\right.
\end{equation}
where $t_0 \ge 0$, $\alpha, \beta : [t_0, +\infty) \to (0, +\infty)$ are continuous, nondecreasing smooth functions, and $g : \mathbb{R}^n \to \mathbb{R}$ is a Fr\'echet differentiable function whose gradient $\nabla g$ is $L_g$-Lipschitz continuous. The initial position and velocity are given by $u_0, v_0 \in \mathbb{R}^n$. This dynamical system is designed to solve the smooth, possibly nonconvex optimization problem
\begin{equation}\label{P}\tag{P}
\inf_{x \in \mathbb{R}^n} g(x).
\end{equation}

Under suitable compatibility and growth conditions on the coefficients $\alpha$ and $\beta$, we establish the global convergence of the trajectories generated by \eqref{dysy} toward critical points of $g$. A distinguishing feature of our framework is that the damping coefficient $\alpha(t)$ is allowed to grow without bound as $t \to +\infty$. This unbounded damping mechanism enables us to surpass the classical exponential convergence rates associated with constant or bounded damping, yielding superlinear-exponent convergence rates for both the trajectories and the objective values.
\subsection{Related results}
A particular instance of \eqref{dysy} was recently examined in \cite{L}, where the system
\begin{equation}\label{eq:intro_laszlo}
\left\{
\begin{aligned}
&\frac{\lambda}{t^2}\ddot{x}(t) + \frac{\gamma}{t}\dot{x}(t) + \nabla g(x(t)) = 0,\\
&x(t_0) = u_0, \quad \dot{x}(t_0) = v_0
\end{aligned}
\right.
\end{equation}
was studied in connection with the minimization problem \eqref{P}. Indeed, dividing the governing equation of \eqref{eq:intro_laszlo} by $\lambda/t^2$ yields
\[
\ddot{x}(t) + \frac{\gamma}{\lambda}t\dot{x}(t) + \frac{t^2}{\lambda}\nabla g(x(t)) = 0,
\]
which corresponds to \eqref{dysy} with $\alpha(t) = \frac{\gamma}{\lambda}t$ and $\beta(t) = \frac{t^2}{\lambda}$. More generally, any second-order system of the form
\[
\lambda(t)\ddot{x}(t) + \gamma(t)\dot{x}(t) + \nabla g(x(t)) = 0
\]
with $\lambda(t) > 0$ can be recast in the form \eqref{dysy} by setting $\alpha(t) = \gamma(t)/\lambda(t)$ and $\beta(t) = 1/\lambda(t)$. In particular, vanishing inertial and damping coefficients of the order $\lambda(t) = \mathcal{O}(t^{-2})$ and $\gamma(t) = \mathcal{O}(t^{-1})$ {correspondingly}, fall directly within the unbounded-damping regime analyzed in this work, provided the induced coefficients satisfy our structural hypotheses.

It is also instructive to compare \eqref{dysy} with the classical heavy-ball system with constant damping studied by B\'egout, Bolte, and Jendoubi \cite{Begout2015}:
\begin{equation}\label{BBJ}
\ddot{y}(s) + \gamma\dot{y}(s) + \nabla g(y(s)) = 0, \qquad y(s_0) = u_0, \quad \dot{y}(s_0) = v_0, \quad \gamma > 0.
\end{equation}
Applying the time reparametrization $s = t^2$ with $x(t) = y(t^2)$ transforms \eqref{BBJ} into
\[
\ddot{x}(t) + \left(2\g t - \frac{1}{t}\right)\dot{x}(t) + 4t^2\nabla g(x(t)) = 0,
\]
which highlights the structural link between constant-damping systems and the scaling $\alpha(t) = \Theta(t)$, $\beta(t) = \Theta(t^2)$ in \eqref{dysy}. In \cite{Begout2015} (see also \cite{BCL}), convergence to a critical point is proved under higher regularity assumptions on $g$ together with a Kurdyka--\L{}ojasiewicz (KL) inequality whose desingularizing function $\varphi$ satisfies $\varphi(s) \ge c\sqrt{s}$ near the origin. By contrast, our analysis requires only that $\nabla g$ is Lipschitz continuous and that the regularized energy
\[
H : \mathbb{R}^n \times \mathbb{R}^n \to \mathbb{R}, \qquad H(u, v) = g(u) + \frac{1}{2}\|u - v\|^2,
\]
satisfies the KL property.

A further comparison can be drawn with the variable-damping dynamical system investigated in \cite{BCL1},
\begin{equation}\label{BCL}
\ddot{x}(t) + \left(\frac{\alpha}{t} + \gamma\right)\dot{x}(t) + \nabla g(x(t)) = 0, \qquad x(t_0) = u_0, \quad \dot{x}(t_0) = v_0,
\end{equation}
where $\alpha, \gamma > 0$. When $\gamma = 0$, \eqref{BCL} reduces to the continuous-time model introduced by Su, Boyd, and Cand\`es \cite{SBC} for Nesterov's accelerated gradient method \cite{Nest1} in the convex setting. Under the time rescaling $t = s^2$, setting $y(s) = x(s^2)$ and $s_0 = \sqrt{t_0}$ turns \eqref{BCL} into
\[
\frac{1}{4s^2}\ddot{y}(s) + \left(\frac{\gamma}{2s} - \frac{2\alpha - 1}{4s^3}\right)\dot{y}(s) + \nabla g(y(s)) = 0, \qquad y(s_0) = u_0, \quad \dot{y}(s_0) = 2s_0 v_0.
\]
In \cite{BCL1}, trajectory convergence is established when the regularized energy $H$ is a KL function. However, even for strongly convex objective functions, the systems in \cite{Begout2015, BCL1} guarantee at best exponential convergence of the form $\|x(t) - \bar{x}\| = \mathcal{O}(e^{-At})$ for some $A > 0$, {as $t \to +\infty$, where $\bar{x}$ is a critical point of $g$}. In contrast, the unbounded damping framework developed here yields convergence rates governed by the integral of $\alpha$, namely
\[
\|x(t) - \bar{x}\| = \mathcal{O}\!\left(e^{-c\int_{t_4}^t \alpha(s)\,ds}\right) \qquad \text{as } t \to +\infty,
\]
{for an appropriate constant $c$ and a starting point $t_4$, which will be defined later in Theorem \ref{thmrates}}, with a corresponding squared-rate decay for the objective residual $g(x(t)) - g(\bar{x})$. In particular, {for positive constants $\alpha_0$ and $A$,} linear damping $\alpha(t) = \alpha_0 t$ yields a Gaussian rate $\mathcal{O}(e^{-At^2})$, polynomial damping $\alpha(t) = \alpha_0 t^p$ ($p > 0$) yields $\mathcal{O}(e^{-At^{p+1}})$, and exponential damping $\alpha(t) = \alpha_0 e^{\lambda t}$ ($\lambda > 0$) achieves a double-exponential rate $\mathcal{O}(e^{-Ae^{\lambda t}})$, {as $t \to +\infty$.}

\subsection{A model result}

As a model result we consider the case when our objective function $g$ is strongly convex, further, we consider the following four classes of second-order dynamical systems, obtained by selecting different profiles for the
time-dependent damping and gradient coefficients.

\begin{enumerate}

\item \textbf{Polynomial damping.}
For $p>0$, $t_0>0$, and $\alpha_0,\beta_0>0$, consider
\begin{equation}\label{DS-polynomial}
\left\{
\begin{aligned}
&\ddot{x}(t)+\alpha_0t^p\dot{x}(t)
+\beta_0t^{2p}\nabla g(x(t))=0,\\
&x(t_0)=u_0,\qquad \dot{x}(t_0)=v_0,
\end{aligned}
\right.
\end{equation}
where
\(
0<\beta_0<
\frac{\sqrt{3}-1}{L_g}\alpha_0^2.
\)

\item \textbf{Exponential damping.}
For $\lambda>0$, $t_0\geq0$, and $\alpha_0,\beta_0>0$,
consider
\begin{equation}\label{DS-exponential}
\left\{
\begin{aligned}
&\ddot{x}(t)+\alpha_0e^{\lambda t}\dot{x}(t)
+\beta_0e^{2\lambda t}\nabla g(x(t))=0,\\
&x(t_0)=u_0,\qquad \dot{x}(t_0)=v_0,
\end{aligned}
\right.
\end{equation}
where
\(
0<\beta_0<
\frac{\sqrt{3}-1}{L_g}\alpha_0^2.
\)

\item \textbf{Logarithmic damping.}
For $\alpha_0,\delta>0$ and
\(
0<\kappa<
\frac{\sqrt{1+2\delta}-1}{L_g\delta^2},
\)
consider
\begin{equation}\label{DS-logarithmic}
\left\{
\begin{aligned}
&\ddot{x}(t)+\alpha_0\log(e+t)\dot{x}(t)+\kappa\alpha_0^2\log^2(e+t)\nabla g(x(t))=0,\\
&x(t_0)=u_0,\qquad \dot{x}(t_0)=v_0.
\end{aligned}
\right.
\end{equation}

\item \textbf{Mixed polynomial--logarithmic--exponential damping.}
For $p,r,\lambda\geq0$, with at least one of these
parameters positive, consider
\begin{equation}\label{DS-mixed}
\left\{
\begin{aligned}
&\ddot{x}(t)+\alpha_0(1+t)^p\log^r(e+t)e^{\lambda t}\dot{x}(t)
+\kappa\alpha_0^2(1+t)^{2p}
\log^{2r}(e+t)e^{2\lambda t}\nabla g(x(t))=0,\\
&x(t_0)=u_0,\qquad \dot{x}(t_0)=v_0,
\end{aligned}
\right.
\end{equation}
where $\alpha_0,\delta>0$ and
\(
0<\kappa<
\frac{\sqrt{1+2\delta}-1}{L_g\delta^2}.
\)
\end{enumerate}
In case the objective function $g$ is strongly convex, we are able to obtain superlinear convergence rates for the trajectories generated by these dynamical systems, and also superlinear rates for the function values in these trajectories.
More precisely the following result hold.
\begin{theorem}\label{thm:strong-convex-superlinear}

Let $g:\mathbb{R}^n\to\mathbb{R}$ be a continuously
differentiable, $\mu$-strongly convex function with
$L_g$-Lipschitz continuous gradient, where $\mu>0$.
Let $u_0,v_0\in\mathbb{R}^n$, and denote the unique minimizer of $g$ by  $x^*.$

Consider the dynamical systems
\eqref{DS-polynomial}--\eqref{DS-mixed}, with the
corresponding parameter restrictions stated above.

Then, for each of these systems, there exist constants
$a,c,T>0$ such that its unique strong global solution
satisfies
\[
\boxed{
\|x(t)-x^*\|
\leq ae^{-c\int_T^t\alpha(s)\,ds},\mbox{ and }g(x(t))-\min g\le  a^2\frac{L_g}{2}e^{-2c\int_T^t\alpha(s)\,ds},
\qquad t\geq T.
}
\]

In particular, the following convergence rates hold:

\begin{enumerate}

\item \textbf{Polynomial damping.}
For $\alpha(t)=\alpha_0t^p$, $p>0$,
\[
\boxed{
\|x(t)-x^*\|
\leq ae^{-ct^{p+1}}\mbox{ and }g(x(t))-\min g\le a^2\frac{L_g}{2}e^{-2ct^{p+1}} .
}
\]
In particular, for $p=1$ we obtain the Gaussian
convergence rate
\[
\|x(t)-x^*\|=\mathcal{O}(e^{-ct^2}) \mbox{ and }g(x(t))-\min g=\mathcal{O}(e^{-2ct^2}).
\]

\item \textbf{Exponential damping.}
For $\alpha(t)=\alpha_0e^{\lambda t}$, $\lambda>0$,
\[
\boxed{
\|x(t)-x^*\|
\leq ae^{-ce^{\lambda t}}\mbox{ and }g(x(t))-\min g\le a^2\frac{L_g}{2}e^{-2ce^{\lambda t}}.
}
\]
Hence, the trajectories converge at a double-exponential
rate.

\item \textbf{Logarithmic damping.}
For $\alpha(t)=\alpha_0\log(e+t)$,
\[
\boxed{
\|x(t)-x^*\|
\leq ae^{-ct\log t}\mbox{ and }g(x(t))-\min g\le a^2\frac{L_g}{2}e^{-2ct\log t}.
}
\]

\item \textbf{Mixed damping.}
For
\[
\alpha(t)=\alpha_0(1+t)^p
\log^r(e+t)e^{\lambda t},
\]
we have
\[
\boxed{
\|x(t)-x^*\|=
\begin{cases}
\mathcal{O}\left(
e^{-ct^{p+1}(\log t)^r}
\right),&\lambda=0,\\[2mm]
\mathcal{O}\left(
e^{-ce^{\lambda t}t^p(\log t)^r}
\right),&\lambda>0
\end{cases}\mbox{ and }g(x(t))-\min g=\begin{cases}
\mathcal{O}\left(
e^{-2ct^{p+1}(\log t)^r}
\right),&\lambda=0,\\[2mm]
\mathcal{O}\left(
e^{-2ce^{\lambda t}t^p(\log t)^r}
\right),&\lambda>0.
\end{cases}
}
\]

\end{enumerate}
\end{theorem}

Consequently, all four classes of dynamical systems
generate trajectories converging to the unique minimizer
of $g$ at superlinear-exponent rates.

\subsection{The outline of the paper}
The remainder of the paper is organized as follows. Section~2 establishes the existence and uniqueness of global strong solutions to \eqref{dysy} and derives an almost-everywhere upper bound on the third-order derivative. Section~3 develops Lyapunov energy estimates and proves that every cluster point of a bounded trajectory is a critical point of $g$. Section~4 establishes full trajectory convergence under the KL property of $H$, with applications to semialgebraic and coercive objectives. Section~5 derives general convergence rates in terms of the \L{}ojasiewicz exponent and identifies explicit parameter regimes yielding superlinear-exponent rates for both trajectories and objective values. Section~6 specializes these results to strongly convex objective functions. Finally, Appendices~A and~B collect auxiliary results on local absolute continuity and the KL property.

\section{Existence, uniqueness and the third order derivative}

We now specify the type of solutions considered in the analysis of the dynamical system \eqref{dysy}.

\begin{definition}\label{str-sol}\rm We say that $x:[t_0,+\infty)\rightarrow \R^n$ is a strong global solution of \eqref{dysy}, if the
following properties are satisfied:

(i) $x,\dot x,\ddot x:[t_0,+\infty)\rightarrow \R^n$ are locally absolutely continuous, in other words, absolutely continuous on each interval $[0,b]$ for $0<b<+\infty$;

(ii) $\ddot{x}(t)+\a(t)\dot{x}(t)+\b(t)\nabla g(x(t))=0$ for almost every $t\geq t_0$;

(iii) $x(t_0)=u_0$ and $\dot{x}(t_0)=v_0$.
\end{definition}

For proving existence and uniqueness of the strong global solutions of \eqref{dysy}, we use the Cauchy-Lipschitz-Picard Theorem for absolutely continuous trajectories (see, for example,
\cite[Proposition 6.2.1]{haraux}, \cite[Theorem 54]{sontag}). The key argument is that one can rewrite \eqref{dysy} as a particular first order dynamical system in a suitably chosen product space (see also \cite{alv-att-bolte-red}).

\begin{theorem}\label{uniq} Assume that $\a,\b$ are continuous functions on $[t_0,+\infty[.$ Then, for each initial pair $u_0,v_0\in \R^n$ there exists a unique strong global solution $x$ of the dynamical system \eqref{dysy}.
\end{theorem}
\begin{proof}

By making use of the notation $X(t)=(x(t),\dot{x}(t))$ the system \eqref{dysy} can be rewritten as a first order dynamical system.
\begin{equation}\label{dysy4}
\left\{
\begin{array}{ll}
\dot{X}(t)=F(t,X(t))\\
X(s_0)=(u_0,v_0),
\end{array}
\right.
\end{equation}
where $F:[t_0,+\infty)\times \R^n\times\R^n\To \R^n\times\R^n,\, F(t,u,v)=\left(v, -\a(t)v-\b(t)\nabla g(u)\big)\right).$
First we show that $F(t,\cdot,\cdot)$ is Lipschitz continuous with a Lipschitz constant $L(t)\in L^1_{loc}([t_0,+\infty))$, for every $t\geq t_0$.

Indeed,
$$\|F(t,u,v)-F(t,\ol u,\ol v)\|=\sqrt{\|v-\ol v\|^2+\left\|\a(t)(\ol v-v)+\b(t)(\nabla g(\ol u)-\nabla g(u))\right\|^2}\le$$
$$\sqrt{\left(1+2\a^2(t)\right)\|v-\ol v\|^2+2L_g^2\b^2(t)\|u-\ol u\|^2}\le$$
$$\sqrt{1+2\a^2(t)+2\b^2(t)}\sqrt{\|v-\ol v\|^2+\|u-\ol u\|^2}=$$
$$\sqrt{1+2\a^2(t)+2\b^2(t)}\|(u,v)-(\ol u,\ol v)\|.$$

Obviously, the Lipschitz constant $L(t)=\sqrt{1+2\a^2(t)+2\b^2(t)}$ is continuous, hence integrable on $[t_0,b]$ for all $0<b<+\infty$, consequently $L(t)\in L^1_{loc}([t_0,+\infty))$.

Next we show that $F(\cdot,u,v)\in L^1_{loc}([t_0,+\infty),{\R^n}\times {\R^n})$ for all $u,v\in {\R^n}.$ Indeed, for $0<b<+\infty$ and $u,v\in {\R^n}$ fixed one has
$$\int_{t_0}^b\|F(t,u,v)\|ds=\int_{t_0}^b \sqrt{\|v\|^2+\left\|\a(t)v+\b(t)\nabla g(u)\right\|^2}dt\le$$
$$\int_{t_0}^b \sqrt{\left(1+2\a^2(t)\right)\|v\|^2+2\b^2(t)\|\nabla g(u)\|^2}dt\le$$
$$\sqrt{\|v\|^2+\|\nabla g(u)\|^2}\int_{t_0}^b\sqrt{1+2\a^2(t)+2\b^2(t)}dt,$$
and the conclusion follows by the continuity of $\sqrt{1+2\a^2(t)+2\b^2(t)}.$

The conclusion of the theorem follows by applying the Cauchy-Lipschitz-Picard theorem to the first order dynamical system \eqref{dysy4}.
\end{proof}

Note that Theorem \ref{uniq} assures that $x^{(3)}$ exists almost everywhere. Let us assume now that $\a,\,\b$ satisfy the following condition. There exist $K_1,K_2>0$ such that for all sufficiently large $t$ it holds
\[\noindent {\bf (L1)}\,\,\left\{\begin{array}{ll}\dot{\a}(t)+\frac{\a(t)\dot{\b}(t)}{{\b}(t)}\le K_1\b(t)\\
\frac{\dot{\b}(t)}{\b(t)}\le K_2\a(t).\end{array}\right.\]

We have the following result.

\begin{lemma}\label{thirdderiv11} Consider $u_0,v_0\in\R^n$ and let $x$ be the unique strong global solution of \eqref{dysy}.

Then  for almost every $t\in[t_0,+\infty)$ one has:
$$\|x^{(3)}(t)\|\le \left( \dot \alpha(t) + \frac{\alpha(t) \dot \beta(t)}{\beta(t)} + L_g \beta(t) \right) \| \dot x(t) \| + \left( \alpha(t) + \frac{\dot \beta(t)}{\beta(t)} \right) \| \ddot x(t) \|.$$
Assume that ${\bf (L1)}$ holds. Then, there exists $C_1,C_2>0$ such that  for almost every $t\in[t_0,+\infty)$ one has:
\begin{align}\label{third-deriv11}
\|x^{(3)}(t)\|\le C_1\b(t) \| \dot x(t) \| + C_2\a(t) \| \ddot x(t) \|.
\end{align}
\end{lemma}

\begin{proof}

From the dynamical system \eqref{dysy} we have for almost every $t\ge t_0$ that

\begin{align*}
    \| \ddot{x}(t + h) - \ddot{x}(t) \| \ = &\left\|  -\alpha(t + h) \dot x(t + h) + \alpha(t) \dot x(t) - \beta(t + h) \nabla g(x(t + h)) + \beta(t) \nabla g(x(t))  \right\| \\
    \leq \ & \| \alpha(t + h) \dot x(t + h) - \alpha(t) \dot x(t) \| + \| (\beta(t + h) - \beta(t)) \nabla g(x(t + h)) \| \\
    &+ \ \beta(t) \| \nabla g(x(t + h)) - \nabla g(x(t)) \| \\
    \leq \ & \| \alpha(t + h) \dot x(t + h) - \alpha(t) \dot x(t) \| \\
    &+ \ \| (\beta(t + h) - \beta(t)) \nabla g(x(t + h)) \| + L_g \beta(t) \| x(t + h) - x(t) \|.
\end{align*}

Hence, for $h>0$ one has

\begin{align*}
   \left \|\frac{ \ddot{x}(t + h) - \ddot{x}(t)}{h}\right \|  \leq \ &\left  \|\frac{ \alpha(t + h) \dot x(t + h) - \alpha(t) \dot x(t)}{h}\right \| \\
    &+ \left \| \frac{\beta(t + h) - \beta(t)}{h} \nabla g(x(t + h)) \right\| + L_g \beta(t)\left \|\frac{ x(t + h) - x(t)}{h}\right \|.
\end{align*}

Thus, by taking the limit as $h \to 0$ in the both sides of the inequality above, we deduce, for almost every $t\ge t_0$ that
\begin{align*}
\| x^{(3)}(t) \| & \leq  \| \dot \alpha(t) \dot x(t) + \alpha(t) \ddot x(t) \| + \left\| \dot{\beta(t)}\n g(x(t)) \right\| + L_g \beta(t) \| \dot x(t)\| \\
&= \| \dot \alpha(t) \dot x(t) + \alpha(t) \ddot x(t) \| + \left\| \dot{\beta(t)}\left(\frac{1}{\b(t)}\ddot{x}(t)+\frac{\a(t)}{\b(t)}\dot{x}(t)\right) \right\| + L_g \beta(t) \| \dot x(t)\|\\
    & \leq \ \left( \dot \alpha(t) + \frac{\alpha(t) \dot \beta(t)}{\beta(t)} + L_g \beta(t) \right) \| \dot x(t) \| + \left( \alpha(t) + \frac{\dot \beta(t)}{\beta(t)} \right) \| \ddot x(t) \|.
    \end{align*}

Now, using ${\bf (L1)}$ we get
$$\| x^{(3)}(t) \| \le C_1\b(t) \| \dot x(t) \| +C_2\a(t)\| \ddot x(t) \|,$$
for some $C_1,C_2>0.$
\end{proof}

\begin{remark}
When \( \alpha(t) = \frac{\gamma}{\lambda} t \) and \( \beta(t) = \frac{t^2}{\lambda} \) we reobtain similar result from \cite{L}, that is,
\[
\| x^{(3)}(t) \| \ \leq \ \left( 3 \frac{\gamma}{\lambda} + \frac{L_g t^2}{\lambda} \right) \| \dot x(t) \| + \left( \frac{\gamma}{\lambda} t + \frac{2}{t} \right) \| \ddot x(t) \|.
\]
\end{remark}

\section{Lyapunov analysis}

In this section we study the asymptotical behaviour of the trajectories of the dynamical system \eqref{dysy}. We recall two technical lemmas which will play an essential role in our analysis, (see  \cite{abbas-att-sv}  Lemma 5.1 and Lemma 5.2, respectively).

\begin{lemma}\label{fejer-cont1} Suppose that $F:[0,+\infty)\rightarrow\R$ is locally absolutely continuous and bounded below and that
there exists $G\in L^1([0,+\infty))$ such that for almost every $t \in [0,+\infty)$ $$\frac{d}{dt}F(t)\leq G(t).$$
Then there exists $\lim_{t\rightarrow \infty} F(t)\in\R$.
\end{lemma}

\begin{lemma}\label{fejer-cont2}  If $1 \leq p < \infty$, $1 \leq r \leq \infty$, $F:[0,+\infty)\rightarrow[0,+\infty)$ is
locally absolutely continuous, $F\in L^p([0,+\infty))$, $G:[0,+\infty)\rightarrow\R$, $G\in  L^r([0,+\infty))$ and
for almost every $t \in [0,+\infty)$ $$\frac{d}{dt}F(t)\leq G(t),$$ then $\lim_{t\rightarrow +\infty} F(t)=0$.
\end{lemma}

 %We start with a very deep abstract result concerning the dynamical system \eqref{dysy}.

{The first result of this section is an auxiliary one, providing important estimates for our future needs.}

\begin{theorem}\label{nocond} Let $x$ be a trajectory generated by \eqref{dysy}. Then,  the function $t\To g\left( x(t) \right)+\frac{1}{2\beta(t)} \|\dot{x}(t)\|^2$ is nonincreasing, the limit
 \begin{equation*}\label{minyag0}
\lim_{t \to +\infty} \left( g\left( x(t) \right) + \frac{1}{2 \beta(t)} \| \dot x(t) \|^2 \right)
\end{equation*}
exists and is finite and furthermore
\begin{equation*}\label{fos}
\sqrt{\frac{\dot{\b}(t)}{2\b^2(t)}+\frac{\a(t)}{\b(t)}}\dot{x}(t),\,\sqrt{\frac{\a(t)}{\b(t)}} \dot x(t),\,\frac{\sqrt{\dot{\b}(t)}}{\b(t)} \dot x(t) \in L^2\left( [t_0, +\infty], \mathbb{R}^n\right).
\end{equation*}
\end{theorem}

\begin{proof} Using the governing differential equation in \eqref{dysy} and multiplying it by $\dot{x}(t)$ we get
$$\frac{1}{\b(t)}\<\ddot{x}(t),\dot{x}(t)\>+\frac{\a(t)}{\b(t)}\|\dot{x}(t)\|^2+\<\n g(x(t)),\dot{x}(t)\>=0.$$
In other words
\begin{equation}\label{fosika}
\frac{d}{dt}\left( g(x(t))+\frac{1}{2\b(t)}\|\dot{x}(t)\|^2\right)+\left(\frac{\dot{\b}(t)}{2\b^2(t)}+\frac{\a(t)}{\b(t)}\right)\|\dot{x}(t)\|^2=0.
\end{equation}
Since $\a,\b$  are positive and nondecreasing we get $\dot{\b}(t)\ge 0$ hence $\frac{d}{dt}\left( g(x(t))+\frac{1}{2\b(t)}\|\dot{x}(t)\|^2\right)\le 0$ which shows
that the function $t\To g\left( x(t) \right)+\frac{1}{2\beta(t)} \|\dot{x}(t)\|^2$ is nonincreasing. Combining the latter with the fact that $g$ is bounded from below, we get that the limit
$$\lim_{t \to +\infty} \left( g\left( x(t) \right) + \frac{1}{2 \beta(t)} \| \dot x(t) \|^2 \right)$$
exists and is finite.

Now, we integrate \eqref{fosika} on an interval $[t_0,T]$ and we take $T\to +\infty$ in order to get
$$\int_{t_0}^{+\infty} \left(\frac{\dot{\b}(t)}{2\b^2(t)}+\frac{\a(t)}{\b(t)}\right)\|\dot{x}(t)\|^2 dt<+\infty.$$
In other words
$$\sqrt{\frac{\dot{\b}(t)}{2\b^2(t)}+\frac{\a(t)}{\b(t)}}\dot{x}(t)\,\sqrt{\frac{\a(t)}{\b(t)}} \dot x(t),\,\frac{\sqrt{\dot{\b}(t)}}{\b(t)} \dot x(t) \in L^2\left( [t_0, +\infty], \mathbb{R}^n\right).$$
\end{proof}

Recall that the set of limit points of the trajectory $x,$ which we denote by  $\omega(x)$, is defined as
$$\omega(x):=\{\ol{x}\in\R^n:\exists t_k\To\infty\mbox{ such that }x(t_k)\To\ol{x},\,k\To + \infty\}.$$
In order to show that every limit point of a  trajectory $x$ is also a critical point of the objective function $g$, we will need the following
$$\lim_{t\to+\infty}\frac{1}{\b(t)}\|\ddot{x}(t)\|=0\mbox{ and }\lim_{t\to+\infty}\frac{\a(t)}{\b(t)}\|\dot{x}(t)\|=0.$$
Indeed, in that case from \eqref{dysy} we get
$$\lim_{t\to+\infty}\n g(x(t))=0.$$

So, we consider an auxiliary function $\g:[t_0,+\infty)\to(0,\infty)$ that satisfies the following conditions:
%$${\bf (L2)}\,\,\,\,\,-\left(1+\dot{\g}(t)\right)\frac{\a(t)}{\b(t)}+\frac{L_g^2\g^3(t)\b(t)}{2}+ L_g \gamma(t) |\dot \gamma(t) + %1|+\left(\frac{\alpha(t) \gamma(t) + \dot \gamma(t) + 1}{2\beta(t)}\right)'\le C_0\left(\frac{\dot{\b}(t)}{2\b^2(t)}+\frac{\a(t)}{\b(t)}\right),$$
%$\mbox{for some }C_0>0\mbox{ and }t\mbox{ big enough.}$
\[{\bf (L2)}\,\,\,\,\,-\left(1+\dot{\g}(t)\right)\frac{\a(t)}{\b(t)}+\frac{L_g^2\g^3(t)\b(t)}{2}+ L_g \gamma(t) |\dot \gamma(t) + 1|+\left(\frac{\alpha(t) \gamma(t) + \dot \gamma(t) + 1}{2\beta(t)}\right)'\le 0 \mbox{ for }t\mbox{ big enough.}\]

\[{\bf (L3)}\,\,\,\,\,\alpha(t) \gamma(t) + \dot \gamma(t) + 1\ge 0\mbox{ for }t\mbox{ big enough.}\]

\[{\bf (L4)}\,\,\,\,\,\frac{1}{\g(t)}\le C_1\a(t)+\left(\frac{C_1}{2}+1\right)\frac{\dot \b(t)}{ \beta(t)}\mbox{ for some }C_1>0\mbox{ and }t\mbox{ big enough.}\]

\begin{theorem}\label{abstr}
Assume that $g$ is bounded from below and for $u_0,v_0\in\R^n$ let $x$ be the unique strong global solution of  the dynamical system \eqref{dysy}. Assume that ${\bf (L2)}-{\bf (L4)}$ hold.

Then

\begin{itemize}
\item[(i)] $\sqrt{\frac{\g(t)}{\b(t)}}  \ddot x(t) \in L^2\left( [t_0, +\infty], \mathbb{R}^n \right)$;
\item[(ii)] there exists the limit $\lim_{t \to +\infty} \left( g\left( \gamma(t) \dot x(t) + x(t) \right) + \frac{\alpha(t) \gamma(t) + \dot \gamma(t) + 1}{2 \beta(t)} \| \dot x(t) \|^2 \right)\in\R$;
\item[(iii)] there exists and is finite the limit $\lim_{t\To+\infty}g(x(t))$;
%\item[(iii)] $\lim_{t\To+\infty}\frac{1}{\b(t)}\ddot{x}(t)=0,$ and $\lim_{t\To+\infty}\frac{\a(t)}{\b(t)}\dot{x}(t)=0$;
%\item[(iv)] $\omega(x)\subseteq \crit(g).$
\end{itemize}
\end{theorem}

\begin{proof}

%For an appropriate choice of $\gamma(t)=\Gamma\sqrt[3]{\frac{\a(t)}{\b^2(t)}},\,\Gamma\in\left[0,\sqrt[3]{\frac{4}{L_g^2}}\right)$ we
We deduce using the $L_g$-Lipschitz continuity of $\nabla g$:

\begin{align}\label{f1}
    &\frac{d}{dt} g\left( \gamma(t) \dot x(t) + x(t) \right)= \langle \nabla g\left( \gamma(t) \dot x(t) + x(t) \right), (\dot \gamma(t) + 1) \dot x(t) + \gamma(t) \ddot x(t) \rangle \\
    \nonumber&= \langle \nabla g\left( \gamma(t) \dot x(t) + x(t) \right) - \nabla g(x(t)), (\dot \gamma(t) + 1) \dot x(t) + \gamma(t) \ddot x(t) \rangle - \frac{1}{\beta(t)} \left\langle \ddot{x}(t) + \alpha(t)\dot{x}(t), (\dot \gamma(t) + 1) \dot x(t) + \gamma(t) \ddot x(t) \right\rangle \\
 \nonumber &\leq \| \nabla g\left( \gamma(t) \dot x(t) + x(t) \right) - \nabla g(x(t))\|\| (\dot \gamma(t) + 1) \dot x(t) + \gamma(t) \ddot x(t)\|  - \frac{\gamma(t)}{\beta(t)} \| \ddot x(t) \|^2 - \frac{\alpha(t) (\dot \gamma(t) + 1)}{\beta(t)} \| \dot x(t) \|^2 \\
\nonumber&- \frac{\alpha(t) \gamma(t) + \dot \gamma(t) + 1}{\beta(t)}\<\dot{x}(t),\ddot{x}(t)\>\\
   \nonumber& \leq L_g \gamma(t) |\dot \gamma(t) + 1| \| \dot x(t) \|^2 + L_g \gamma^2(t) \| \dot x(t)\|\| \ddot x(t) \|- \frac{\gamma(t)}{\beta(t)} \| \ddot x(t) \|^2 - \frac{\alpha(t) (\dot \gamma(t) + 1)}{\beta(t)} \| \dot x(t) \|^2 \\
   \nonumber &-  \frac{\alpha(t) \gamma(t) + \dot \gamma(t) + 1}{\beta(t)} \langle \dot x(t), \ddot x(t) \rangle.
   \end{align}
Now, observe that
\begin{equation}\label{f11}
-  \frac{\alpha(t) \gamma(t) + \dot \gamma(t) + 1}{\beta(t)} \langle \dot x(t), \ddot x(t) \rangle=-\frac{d}{dt}\left( \frac{\alpha(t) \gamma(t) + \dot \gamma(t) + 1}{2\beta(t)} \|\dot{x}(t)\|^2\right)+\left(\frac{\alpha(t) \gamma(t) + \dot \gamma(t) + 1}{2\beta(t)}\right)' \|\dot{x}(t)\|^2
\end{equation}
and for  $a:[t_0,+\infty)\to(0,+\infty),\,a(t)=L_g\g(t)\b(t)$  one has
\begin{align}\label{f12}
 L_g \gamma^2(t) \| \dot x(t)\|\| \ddot x(t) \|&\le\frac{L_g\g^2(t)a(t)}{2}\|\dot{x}(t)\|^2+\frac{L_g\g^2(t)}{2a(t)}\|\ddot{x}(t)\|^2\\
&\nonumber=\frac{L_g^2\g^3(t)\b(t)}{2}\|\dot{x}(t)\|^2+\frac{\g(t)}{2\b(t)}\|\ddot{x}(t)\|^2,
\end{align}
for all  $t\in[t_0,+\infty).$
Hence, by injecting \eqref{f11} and \eqref{f12} in \eqref{f1}, from ${\bf (L2)}$ we get that there exist  $t_1\ge t_0$ such that for all $t\ge t_1$ the following holds

\begin{align}\label{f2}
    &\frac{d}{dt} \left(g\left( \gamma(t) \dot x(t) + x(t) \right)+\frac{\alpha(t) \gamma(t) + \dot \gamma(t) + 1}{2\beta(t)} \|\dot{x}(t)\|^2\right)\le-\frac{\g(t)}{2\b(t)}\|\ddot{x}(t)\|^2+\\
    \nonumber&\left(-\left(1+\dot{\g}(t)\right)\frac{\a(t)}{\b(t)}+\frac{L_g^2\g^3(t)\b(t)}{2}+ L_g \gamma(t) |\dot \gamma(t) + 1|+\left(\frac{\alpha(t) \gamma(t) + \dot \gamma(t) + 1}{2\beta(t)}\right)' \right)\|\dot{x}(t)\|^2.
    \end{align}

    In other words, for all $t\ge t_1$, one has
    \begin{align}\label{ineq}
    &\frac{d}{dt} \left(g\left( \gamma(t) \dot x(t) + x(t) \right)+\frac{\alpha(t) \gamma(t) + \dot \gamma(t) + 1}{2\beta(t)} \|\dot{x}(t)\|^2\right)+\frac{\g(t)}{2\b(t)}\|\ddot{x}(t)\|^2+\\
    \nonumber&\left(\left(1+\dot{\g}(t)\right)\frac{\a(t)}{\b(t)}-\frac{L_g^2\g^3(t)\b(t)}{2}- L_g \gamma(t) |\dot \gamma(t) + 1|-\left(\frac{\alpha(t) \gamma(t) + \dot \gamma(t) + 1}{2\beta(t)}\right)' \right)\|\dot{x}(t)\|^2\le 0.
    \end{align}

    By integrating the relation \eqref{ineq} on $[t_1, T]$ we obtain

\begin{align}\label{ineqfos}
    &g\left( \gamma(T) \dot x(T) + x(T) \right) + \frac{\alpha(T) \gamma(T) + \dot \gamma(T) + 1}{2 \beta(T)} \| \dot x(T) \|^2+ \int_{t_1}^T \frac{\g(t)}{2\b(t)}\|\ddot{x}(t)\|^2 dt \\
    \nonumber &+\int_{t_1}^T\left(\left(1+\dot{\g}(t)\right)\frac{\a(t)}{\b(t)}-\frac{L_g^2\g^3(t)\b(t)}{2}- L_g \gamma(t) |\dot \gamma(t) + 1|-\left(\frac{\alpha(t) \gamma(t) + \dot \gamma(t) + 1}{2\beta(t)}\right)' \right)\|\dot{x}(t)\|^2 dt\\
    \nonumber\leq \ &g\left( \gamma(t_1) \dot x(t_1) + x(t_1) \right) + \frac{\alpha(t_1) \gamma(t_1) + \dot \gamma(t_1) + 1}{2 \beta(t_1)} \| \dot x(t_1) \|^2.
\end{align}

Taking into account the boundedness of $g$ from below and conditions {${\bf (L2)}$ and} ${\bf (L3)}$ we deduce by taking the limit as $T \to +\infty$ that

\[
\int_{t_1}^{+\infty}\frac{\g(t)}{2\b(t)}\|\ddot{x}(t)\|^2 dt < +\infty,
\]
which means nothing but
\[
\sqrt{\frac{\g(t)}{\b(t)}}  \ddot x(t) \in L^2\left( [t_1, +\infty], \mathbb{R}^n \right).
\]

Furthermore, according to \eqref{ineq} we have that the function  $t \To g\left( \gamma(t) \dot x(t) + x(t) \right) + \frac{\alpha(t) \gamma(t) + \dot \gamma(t) + 1}{2 \beta(t)} \| \dot x(t) \|^2$
is decreasing and the fact that $g$ is bounded from below {together with ${\bf (L3)}$} ensures that the limit
\begin{equation}\label{minya}
\lim_{t \to +\infty} \left( g\left( \gamma(t) \dot x(t) + x(t) \right) + \frac{\alpha(t) \gamma(t) + \dot \gamma(t) + 1}{2 \beta(t)} \| \dot x(t) \|^2 \right)
\end{equation}
exists and is actually finite. However, we already know from Theorem \ref{nocond} that  the limit
\begin{equation*}
\lim_{t \to +\infty} \left( g\left( x(t) \right) + \frac{1}{2 \beta(t)} \| \dot x(t) \|^2 \right)
\end{equation*}
exists and is finite. Hence, our next goal is to show that
 $\lim_{t \to +\infty} \frac{1}{ \beta(t)} \| \dot x(t) \|^2=0$ and for this purpose we deploy Lemma \ref{fejer-cont2}.

 Consider

\begin{align*}
    \frac{d}{dt} \left( \frac{1}{ \beta(t)} \| \dot x(t) \|^2 \right) =\ &\frac{2}{\beta(t)} \langle \dot x(t), \ddot x(t) \rangle - \frac{\dot{\b}(t)}{ \beta^2(t)} \| \dot x(t) \|^2 \\
    \leq \ & - \frac{\dot \b(t)}{ \beta^2(t)} \| \dot x(t) \|^2 + \frac{1}{\beta(t)}\left(\frac{1}{\g(t)}\|\dot{x}(t)\|^2+\g(t) \| \ddot x(t) \|^2\right) \\
    =\ &\frac{\g(t)}{\b(t)} \| \ddot x(t)\|^2+\left(-\frac{\dot \b(t)}{ \beta^2(t)}+ \frac{1}{\b(t)\g(t)}\right)\| \dot x(t) \|^2.
\end{align*}

According to ${\bf (L4)}$, {$(i)$ of the current Theorem} and Theorem \ref{nocond}, the right-hand side of the previous relation belongs to $L^1\left( [t_1', +\infty], \mathbb{R}\right)$, for some $t_1'\ge t_0$, hence Lemma \ref{fejer-cont2} yields  $\lim_{t \to +\infty} \frac{1}{ \beta(t)} \| \dot x(t) \|^2=0$.
Consequently, the limit
\begin{equation*}
\lim_{t \to +\infty} g\left( x(t) \right)
\end{equation*}
exists and is finite.
\end{proof}

{In order to show that the limit points of the generated trajectories are actually critical points of $g$ (as well as some other important estimates), merely $\lim_{t \to +\infty} \frac{1}{ \beta(t)} \| \dot x(t) \|^2=0$ in not enough. We will to prove next two additional auxiliary facts and two additional assumptions.} In order to show that  $\lim_{t \to +\infty} \frac{\alpha(t)}{ \beta(t)} \| \dot x(t) \|=0$, we assume:
 \[{\bf (L5)}\,\,\,\,\frac{ 2\dot \alpha(t) \a(t)\b(t) -2 \alpha^2(t) \dot \b(t) }{ \beta^2(t)}+\frac{\alpha^4(t)}{\beta^2(t)\g(t)}\le C_2\left(\frac{\dot{\b}(t)}{2\b(t)}+\a(t)\right)\mbox{ for some }C_2>0\mbox{ and }t\mbox{ big enough.}
 \]
 In order to show that $ \lim_{t \to +\infty} \frac{1}{\beta(t)} \| \ddot x(t) \| = 0$ we additionally require:
   \[{\bf (L6)}\,\,\,\,\frac{(\dot{\a}(t)\b(t)+\a(t)\dot{\b}(t)+L_g\b^2(t))^2}{\g(t)\b^4(t)}\le C_3\left(\frac{\dot{\b}(t)}{2\b(t)}+\a(t)\right)\mbox{ for some }C_3>0\mbox{ and }t\mbox{ big enough}.\]
   \[{\bf (L7)}\,\,\,\,C_5\g(t)\le\frac{\a(t)}{\b(t)}\le C_4\g(t)\mbox{ for some }C_4,C_5>0\mbox{ and }t\mbox{ big enough}.\]
  
We have the following result.
  
 \begin{theorem}\label{abstr1}
 Assume that $g$ is bounded from below and for $u_0,v_0\in\R^n$ let $x$ be the unique strong global solution of  the dynamical system \eqref{dysy}. Assume that ${\bf (L2)}-{\bf (L7)}$ hold.

Then

\begin{itemize}
\item[(i)] $\lim_{t\To+\infty}\frac{1}{\b(t)}\ddot{x}(t)=0,$ and $\lim_{t\To+\infty}\frac{\a(t)}{\b(t)}\dot{x}(t)=0$;
\item[(ii)] $\omega(x)\subseteq \crit(g).$
\end{itemize}
\end{theorem}
\begin{proof}

We have
\begin{align*}
    &\frac{d}{dt} \left( \frac{\alpha^2(t)}{ \beta^2(t)} \| \dot x(t) \|^2 \right) =2\frac{\alpha^2(t)}{\beta^2(t)} \langle \dot x(t), \ddot x(t) \rangle + \frac{ 2\dot \alpha(t) \a(t)\b(t) -2 \alpha^2(t) \dot \b(t) }{ \beta^3(t)} \| \dot x(t) \|^2 \\
    \leq \ &  \frac{ 2\dot \alpha(t) \a(t)\b(t) -2 \alpha^2(t) \dot \b(t) }{ \beta^3(t)} \| \dot x(t) \|^2 + \frac{\alpha^2(t)}{\beta^2(t)}\left(\frac{\a^2(t)}{\g(t)\b(t)}\|\dot{x}(t)\|^2+\frac{\g(t)\b(t)}{\a^2(t)} \| \ddot x(t) \|^2\right) \\
    \le\ &{\frac{\g(t)}{\b(t)}} \| \ddot x(t)\|^2+C_2\left(\frac{\dot{\b}(t)}{2\b^2(t)}+\frac{\a(t)}{\b(t)}\right)\| \dot x(t) \|^2,
\end{align*}
where above we used $({\bf L5}).$ The right-hand side of the previous relation belongs to $L^1\left( [t_1, +\infty], \mathbb{R} \right)$ {due to Theorems \ref{nocond} and \ref{abstr}}, hence

\begin{equation}\label{lim_0velocity}
    \lim_{t \to +\infty} \frac{\alpha(t)}{\beta(t)} \| \dot x(t) \| = 0 \quad \text{by Lemma \ref{fejer-cont2}}.
\end{equation}

In order to show that $ \lim_{t \to +\infty} \frac{1}{\beta(t)} \| \ddot x(t) \| = 0$ we use $({\bf L6})$ and $({\bf L7})$ 
    combined with  Lemma \ref{thirdderiv11}.  We get
\begin{align*}
    &\frac{d}{dt} \left( \frac{1}{\beta^2(t)} \| \ddot x(t) \|^2 \right)=-\frac{2\dot{\b}(t)}{\b^3(t)}  \| \ddot x(t) \|^2+2  \frac{1}{\beta^2(t)} \< \ddot x(t),x^{(3)}(t)\>\\
     &\le-\frac{2\dot{\b}(t)}{\b^3(t)}  \| \ddot x(t) \|^2+ \frac{2}{\beta^2(t)}\|\ddot{x}(t)\|\left(\left( \dot \alpha(t) + \frac{\alpha(t) \dot \beta(t)}{\beta(t)} + L_g \beta(t) \right) \| \dot x(t) \| + \left( \alpha(t) + \frac{\dot \beta(t)}{\beta(t)} \right) \| \ddot x(t) \|\right)\\
       &\le  \left( \frac{\dot \alpha(t)}{\b^2(t)} + \frac{\alpha(t) \dot \beta(t)}{\beta^3(t)} + \frac{L_g}{\beta(t)} \right) \left(\frac{\dot{\a}(t)\b(t)+\a(t)\dot{\b}(t)+L_g\b^2(t)}{\g(t)\b^2(t)}\|\dot{x}(t)\|^2+\frac{\g(t)\b^2(t)}{\dot{\a}(t)\b(t)
       +\a(t)\dot{\b}(t)+L_g\b^2(t)}\|\ddot{x}(t)\|^2\right)\\
       &+\frac{2\a(t)}{\beta^2(t)}  \|\ddot{x}(t)\|^2= \frac{(\dot{\a}(t)\b(t)+\a(t)\dot{\b}(t)+L_g\b^2(t))^2}{\g(t)\b^5(t)}\|\dot{x}(t)\|^2+\left(\frac{\g(t)}{\b(t)}
       +\frac{2\a(t)}{\beta^2(t)}\right)\|\ddot{x}(t)\|^2\\
       &\le  C_3\left(\frac{\dot{\b}(t)}{2\b^2(t)}+\frac{\a(t)}{\b(t)}\right)\|\dot{x}(t)\|^2+\left(1+\frac{C_4}{2}\right)\frac{\g(t)}{\b(t)}\|\ddot{x}(t)\|^2.
\end{align*}
Simiraliry, the right-hand side of the previous relation belongs to $L^1\left( [t_1, +\infty], \mathbb{R} \right)$ {due to Theorems \ref{nocond} and \ref{abstr}}, hence

\begin{equation}\label{lim_0acceleration}
    \lim_{t \to +\infty} \frac{1}{\beta(t)} \| \ddot x(t) \| = 0 \quad \text{by Lemma \ref{fejer-cont2}}.
\end{equation}

By combining \eqref{lim_0velocity},  \eqref{lim_0acceleration} and \eqref{dysy} we get
\begin{equation}\label{lim_0gradient}
    \lim_{t \to +\infty}\n g(x(t))=0,\mbox{ hence }\omega(x)\subseteq \crit(g).
\end{equation}

\end{proof}

\begin{remark}\label{abstgamma} Note that under the assumptions of Theorem \ref{abstr1} one has
$$\lim_{t\To+\infty}{\g(t)}\dot{x}(t)=0.$$
This follows directly from condition {\bf (L7)}.
\end{remark}

In what follows we introduce a regularization function of $g$ and show that strong conection between the criticl points of $g$ and the critical points of this  regularization function.
Therefore, consider \( H: \mathbb{R}^n \times \mathbb{R}^n \to \mathbb{R}, \ H(x, y) = g(x) + \frac{1}{2} \| x - y \|^2 \), and denote
    \[
    u(t) \ = \ \gamma(t) \dot x(t) + x(t)
    \]
    and
    \[
    v(t) \ = \ \left( \sqrt{\frac{\alpha(t) \gamma(t) + \dot \gamma(t) + 1}{\beta(t)}} + \gamma(t) \right) \dot x(t) + x(t).
    \]
    Note that $H(u(t),v(t))=g(\gamma(t) \dot x(t) + x(t))+\frac{\alpha(t) \gamma(t) + \dot \gamma(t) + 1}{2\beta(t)} \|\dot x(t)\|^2.$
    
    Assume the following:
    \[{\bf (L8)}\,\,\,\,  \sqrt{\frac{\alpha(t) \gamma(t) + \dot \gamma(t) + 1}{\beta(t)}}\le C_6\frac{\a(t)}{\b(t)}\mbox{ for some }C_6>0\mbox{ and }t\mbox{ big enough.}\]

The next theorem is a consequence of Theorem \ref{nocond}, Theorem \ref{abstr} and Theorem \ref{abstr1}.

\begin{theorem}\label{thm}
Assume that $g$ is bounded from below and for $u_0,v_0\in\R^n$ let $x$ be the unique strong global solution of  the dynamical system \eqref{dysy}. Assume that ${\bf (L2)}-{\bf (L8)}$ hold.

Then we have the following:
    \begin{itemize}

        \item[(i)] \( \omega(u) = \omega(v) = \omega(x) \);

        \item[(ii)] \( \frac{d}{dt} H(u(t), v(t)) \ \leq \ -A(t)\|\dot{x}(t)\|^2-B(t)\|\ddot{x}(t)\|^2 \ \leq \ 0 \) for all \( t \geq t_1 \),
        
        where $A(t)=\left(1+\dot{\g}(t)\right)\frac{\a(t)}{\b(t)}-\frac{L_g^2\g^3(t)\b(t)}{2}- L_g \gamma(t) |\dot \gamma(t) + 1|-\left(\frac{\alpha(t) \gamma(t) + \dot \gamma(t) + 1}{2\beta(t)}\right)'$ and $B(t)=\frac{\g(t)}{2\b(t)}.$

        \item[(iii)] \( \lim_{t \to +\infty} H(u(t), v(t)) \ = \ \lim_{t \to +\infty} g(x(t)) \in \mathbb{R} \);

        \item[(iv)] \( H \) is finite and constant on \( \omega(u, v) \ = \ \left\{ (\bar x, \bar x) \in \mathbb{R}^n \times \mathbb{R}^n: \bar x\in \omega(x) \right\} \);

        \item[(v)] \( \| \nabla H(u(t), v(t)) \| \ \leq \ q(t) \| \ddot x(t) \| + p(t) \| \dot x(t) \| \) for all \( t \geq t_1 \)
        
        where \( q(t) = \frac{1}{\beta(t)} \) and \( p(t) = L_g \gamma(t) + \frac{\alpha(t)}{\beta(t)} + 2 \sqrt{\frac{\alpha(t) \gamma(t) + \dot \gamma(t) + 1}{\beta(t)}} \);

        \item[(vi)] \( \omega(u, v) \subseteq \crit(H) \).

    Assume furthermore that $x$ is bounded. Then

        \item[(vii)] \( \omega(u, v) \) is nonempty and compact

        and

        \item[(viii)] \( \lim_{t \to +\infty} \dist\left( \big( u(t), v(t) \big), \omega(u, v) \right) = 0 \).

    \end{itemize}

\end{theorem}

\begin{proof}

We have by  $({\bf L7}),({\bf L8})$ and Theorem \ref{abstr1} that
\begin{equation}\label{lim}
\lim_{t \to +\infty} \gamma(t) \dot x(t) \ = \ \lim_{t \to +\infty} \sqrt{\frac{\alpha(t) \gamma(t) + \dot \gamma(t) + 1}{\beta(t)}}  \dot x(t)  \ = \ 0,
\end{equation}
hence the first claim follows. The second claim follows from \eqref{f2} and the third one from $({\bf L7})$, $({\bf L8})$ and Theorem \ref{abstr1} correspondingly.

For the proof of the fourth claim  consider \( (\bar x, \bar y) \in \omega(u, v) \). Then, there exists a sequence \( t_k \to +\infty \) as \( k \to +\infty \), such that \( u(t_k) \to \bar x \), as \( k \to +\infty \), and \( v(t_k) \to \bar y \), as \( k \to +\infty \). According to Theorem \ref{abstr1} and taking into account $({\bf L7}),({\bf L8})$ and the forms of \(u\) and \(v\) we deduce that \( \bar x = \bar y \) and \( \bar x \in \omega(x) \). That means \( \omega(u, v) = \left\{ (\bar x, \bar x) \in \mathbb{R}^n \times \mathbb{R}^n: \bar x \in \omega(x) \right\} \). Let\( (\bar x, \bar x) \in \omega(u, v) \). Then there exists a sequence \( t_k \to +\infty \) as \( k \to +\infty \), such that \( u(t_k) \to \bar x \) as \( k \to +\infty \), and \( v(t_k) \to \bar x \) as \( k \to +\infty \). According to \( (iii) \) the limit \( \lim_{t \to +\infty} H(u(t), v(t)) \) exists and is finite, therefore \( H(\bar x, \bar x) = \lim_{k \to +\infty} H(u(t_k), v(t_k)) = \lim_{t \to +\infty} H(u(t), v(t)) \) exists and is finite. That means nothing but that \(H\) is constant on \( \omega(u, v) \).

To prove the fifth claim we use \eqref{dysy} and we get

\begin{align*}
    \| \nabla H(u(t), v(t)) \| \ &=\|(\n g(u(t))+u(t)-v(t),v(t)-u(t))\|\\    
    &\leq \ \| \nabla g(u(t)) \| + 2 \| u(t) - v(t) \| \ \leq \ \| \nabla g(u(t)) - \nabla g(x(t)) \| + \| \nabla g(x(t)) \| + 2 \| u(t) - v(t) \| \\
    &\leq \ L_g \gamma(t) \| \dot x(t) \| + \left\| -\frac{\ddot{x}(t) + \alpha(t)\dot{x}(t)}{\beta(t)} \right\| + 2 \| u(t) - v(t) \| \\
    &\leq \ \frac{1}{\beta(t)} \| \ddot x(t) \| + \left( L_g \gamma(t) + \frac{\alpha(t)}{\beta(t)} + 2 \sqrt{\frac{\alpha(t) \gamma(t) + \dot \gamma(t) + 1}{\beta(t)}} \right) \| \dot x(t) \| \\
    &= \ q(t) \| \ddot x(t) \| + p(t) \| \dot x(t) \|
\end{align*}

The rest of the proof goes in line with \cite{L}. Notice that \( \crit(H) = \left\{ (x, y) \in \mathbb{R}^n \times \mathbb{R}^n: (0, 0) \in \nabla H(x, y) \right\} = \left\{ (\bar x, \bar x) \in \mathbb{R}^n \times \mathbb{R}^n: \bar x \in \crit(g) \right\} \), thus, using the fact established earlier, namely, \\ \( \omega(u, v) = \left\{ (\bar x, \bar x) \in \mathbb{R}^n \times \mathbb{R}^n: \bar x \in \omega(x) \right\} \), and Theorem \ref{thm} item \( (iv) \) one has \( \omega(u, v) \subseteq \crit(H) \), which is the sixth claim.

Claim seven. Since \(x\) is bounded by assumption, by \eqref{lim} we get that \(u\) and \(v\) are also bounded. By invoking Bolzano-Weierstrass Theorem we obtain that the set of limit points \( \omega(u, v) \) is nonempty. Further, \( \omega(u, v) = \left\{ (\bar x, \bar x) \in \mathbb{R}^n \times \mathbb{R}^n: \bar x \in \omega(x) \right\} \) and since \(\omega(x)\) is bounded, it is enough to show that it is closed. In order to do so, we consider a sequence \( \{\bar{x}_n\}_{n \geq 1} \subseteq \omega(x) \) and assume that \( \lim_{n \to +\infty} \bar{x}_n = x^* \). We will show that \( x^* \in \omega(x) \). For every \( n \geq 1 \) there exists the sequence \( t_k^n \to +\infty \), as \( k \to +\infty \), such that \( \lim_{k \to _\infty} x(t_k^n) = \bar{x}_n \). Then, \( \lim_{n \to _\infty} \bar{x}_n = x^* \), therefore for every \( \varepsilon > 0 \) there exists \( N_\varepsilon \in \mathbb{N} \) such that for all \( n \geq N_\varepsilon \) one has \( \| \bar{x}_n - x^* \| \leq \frac{\varepsilon}{2} \). Fixing \( \varepsilon > 0 \) and \( n \geq 1 \), we deduce, by the fact that \( \lim_{k \to _\infty} x(t_k^n) = \bar{x}_n \), the existence of \( k(n, \varepsilon) \in \mathbb{N} \) such that for all \( k \geq k(n, \varepsilon) \) one has \( \| x(t_k^n) - \bar{x}_n \| \leq \frac{\varepsilon}{2} \). Combining the last two inequalities together, we obtain
\[
\| x(t_k^n) - x^* \| \ < \ \varepsilon, \ \forall n \geq N_\varepsilon, \ \forall k \geq k(n, \varepsilon).
\]
Since \( t_k^n \to +\infty \), as \( k \to +\infty \), for all \( n \geq 1 \), we obtain that for every \( n \geq 1 \) there exists \( N_n \in \mathbb{N} \) such that \( t_k^n > n \) for all \( k \geq N_n \). Consider now for every \( n \geq 1 \) the sequence \( s_n = t_{k_n}^n \) and \( k_n > \max\{ n, N_n, k(n, \varepsilon) \} \). It is easy to see that \( s_n \to +\infty \), as \( n \to +\infty \), and \( \| x(s_n) - x^* \| < \varepsilon \ \forall n \geq N_\varepsilon \). Therefore,
\[
\lim_{n \to +\infty} x(s_n) = x^* \text{ or } x^* \in \omega(x).
\]

The last claim follows from the fact that
\[
0 \ \leq \ \lim_{t \to +\infty} \dist((u(t), v(t)), \omega(u, v)) \ \leq \ \lim_{k \to +\infty} \dist((u(t_k), v(t_k)), \omega(u, v)) \ = \ 0.
\]

\end{proof}

\section{Convergence under the KL inequality}

{This section is devoted to proving the convergence of the generated trajectories of \eqref{dysy}} under the assumption of KL property on $H$. 
  Assume additionally the following:
  \[{\bf( L9)}\, 0<\left(1+\dot{\g}(t)-C_7\right)\frac{\a(t)}{\b(t)}-\frac{L_g^2\g^3(t)\b(t)}{2}- L_g \gamma(t) (|\dot \gamma(t) + 1|+C_7)-\left(\frac{\alpha(t) \gamma(t) + \dot \gamma(t) + 1}{2\beta(t)}\right)'-2C_7\sqrt{\frac{\alpha(t) \gamma(t) + \dot \gamma(t) + 1}{\beta(t)}}\]
   and
\[\frac{C_7}{2-4C_8}\frac{\b(t)}{\g(t)}\left(\frac{C_7}{\b(t)}+C_8\g(t)\left(L_g \gamma(t) + \frac{\alpha(t)}{\beta(t)} + 2 \sqrt{\frac{\alpha(t) \gamma(t) + \dot \gamma(t) + 1}{\beta(t)}}\right)\right)^2\le\]
\[\left(1+\dot{\g}(t)-C_7\right)\frac{\a(t)}{\b(t)}-\frac{L_g^2\g^3(t)\b(t)}{2}- L_g \gamma(t) (|\dot \gamma(t) + 1|+C_7)-\left(\frac{\alpha(t) \gamma(t) + \dot \gamma(t) + 1}{2\beta(t)}\right)'-2C_7\sqrt{\frac{\alpha(t) \gamma(t) + \dot \gamma(t) + 1}{\beta(t)}}\]
 $\mbox{ for some }C_7>0,\frac12>C_8>0\mbox{ and }t\mbox{ big enough.}$

We will first prove an auxiliary Lemma.
\begin{lemma}\label{fosikacska} Under the assumptions of Theorem \ref{thm} and $({\bf L9})$,  there exists $a,b>0$ such that for all sufficiently large $t$ one has
$$\frac{A(t) \| \dot x(t) \|^2+B(t) \| \ddot x(t) \|^2}{q(t) \| \ddot x(t) \| + p(t) \| \dot x(t) \|} \ge a\| \dot x(t) \|+b\g(t)\| \ddot x(t) \|.$$
\end{lemma}
\begin{proof}
Note that $\frac{A(t) \| \dot x(t) \|^2+B(t) \| \ddot x(t) \|^2}{q(t) \| \ddot x(t) \| + p(t) \| \dot x(t) \|} \ge a\| \dot x(t) \|+b\g(t)\| \ddot x(t) \|$ can equivalently be written as
$$(A(t)-ap(t))\|\dot{x}(t)\|^2-(aq(t)+b\g(t)p(t))\|\dot{x}(t)\|\|\ddot{x}(t)\|+(B(t)-b\g(t)q(t))\|\ddot{x}(t)\|^2\ge 0,$$
for all sufficiently large $t$.
Hence, if there exists $t'>t_0$ such that 

$0<a<\frac{A(t)}{p(t)}$, $0<b<\frac{B(t)}{\g(t)q(t)}$ and
$(aq(t)+b\g(t)p(t))^2-4(A(t)-ap(t))(B(t)-b\g(t)q(t))\le 0$ for all $t>t'$ then our claim follows. 

Note that  $A(t)=\left(1+\dot{\g}(t)\right)\frac{\a(t)}{\b(t)}-\frac{L_g^2\g^3(t)\b(t)}{2}- L_g \gamma(t) |\dot \gamma(t) + 1|-\left(\frac{\alpha(t) \gamma(t) + \dot \gamma(t) + 1}{2\beta(t)}\right)'$ and $B(t)=\frac{\g(t)}{2\b(t)}$ and \( q(t) = \frac{1}{\beta(t)} \) and \( p(t) = L_g \gamma(t) + \frac{\alpha(t)}{\beta(t)} + 2 \sqrt{\frac{\alpha(t) \gamma(t) + \dot \gamma(t) + 1}{\beta(t)}} \), hence we just use $({\bf L9})$ with $a=C_7$ and $b=C_8.$
\end{proof}

Next Theorem states the main result of the section.

\begin{theorem}\label{convergence}

    Let \(x\) be a unique strong global solution of \eqref{dysy}. Consider the following regularization of \(g\):
    \[
    H: \mathbb{R}^n \times \mathbb{R}^n \to \mathbb{R}, \ H(x, y) = g(x) + \frac{1}{2} \| x - y \|^2
    \]
    and suppose that \(H\) is a KL-function and \(x\) is bounded. Then under the assumptions of Theorem \ref{thm} the following holds
    \begin{itemize}

        \item[(i)] \( \| \dot x(t) \| \in L^1\left( [t_0, +\infty), \mathbb{R} \right) \);

        \item[(ii)] \( \gamma(t) \| \ddot x(t) \| \in L^1\left( [t_0, +\infty), \mathbb{R} \right) \)

        and

        \item[(iii)] there exists \( \bar x \in \crit(g) \) such that \( \lim_{t \to +\infty} x(t) = \bar x \).

    \end{itemize}

\end{theorem}

\begin{proof}

The proof of this result is inspired by \cite{L}. Consider a point \( (\bar x, \bar x) \in \omega(u, v) \) and notice that \( \lim_{t \to +\infty} H(u(t), v(t)) = H(\bar x, \bar x) \) and $\bar x\in\omega(x)\subseteq \crit g$. From now on we must consider two cases.
\vskip0.3cm
\textbf{Case \(1\).} Suppose there exists \( t_2 \geq t_1 \) such that \( H(u(t_2), v(t_2)) = H(\bar x, \bar x) \), where \(t_1\) was defined in Theorem \ref{abstr}. Invoking Theorem \ref{thm}\((ii)\) we obtain
\[
\frac{d}{dt} H(u(t), v(t)) \ \leq \ -A(t)\|\dot{x}(t)\|^2-B(t)\|\ddot{x}(t)\|^2 \ \leq \ 0, \text{ for all } t \geq t_2
\]
and therefore
\[
H(u(t), v(t)) \ \leq \ H(\bar x, \bar x), \text{ for all } t \geq t_2.
\]
But
\[
H(u(t), v(t)) \ \geq \ \lim_{t \to +\infty} H(u(t), v(t)) \ = \ H(\bar x, \bar x), \text{ for all } t \geq t_2,
\]
so
\[
H(u(t), v(t)) \ = \ H(\bar x, \bar x), \text{ for all } t \geq t_2,
\]
which means that \(H\) is constant and \( \frac{d}{dt} H(u(t), v(t)) = 0\), for all \( t \geq t_2 \), which, in turn, leads to
\[
0 \ \leq \ -A(t)\|\dot{x}(t)\|^2-B(t)\|\ddot{x}(t)\|^2 \ \leq \ 0, \text{ for all } t \geq t_2
\]
and since \( A(t), B(t) > 0 \)  for all values of \(t\), we deduce
\[
\| \dot x(t) \| \ = \ \| \ddot x(t) \| \ = \ 0, \text{ for all } t \geq t_2.
\]
Thus, \( \| \dot x(t) \|, \| \ddot x(t) \| \in L^1\left( [t_0, +\infty), \mathbb{R}^n \right) \) and \( x(t) = \bar x \) is constant on \( [t_2, +\infty) \). Finally, it is clear that \( \gamma(t) \| \ddot x(t) \| \in L^1\left( [t_0, +\infty), \mathbb{R}^n \right) \) since \( \gamma \) is continuous on \( [t_0, t_2] \) and $\int_{t_0}^{+\infty} \g(t)\|\ddot{x}(t)\|dt=\int_{t_0}^{t_2} \g(t)\|\ddot{x}(t)\|dt.$
\vskip0.3cm
\textbf{Case \( 2 \).}
Suppose that \( H(u(t), v(t)) > H(\bar x, \bar x) \) for every \( t \geq t_1 \). According to Theorem \ref{thm}\((iv)\) and \((vii)\), \( \omega(u, v) \) is non-empty and compact and \(H\) is constant on \( \omega(u, v) \). Now, since \(H\) is a KL-function, we deduce from Lemma \ref{unif-KL-property}  {(see Appendix B)} that there exist \( \varepsilon, \eta > 0 \) and \( \varphi \in \Theta_\eta \) such that for all \( (\bar x, \bar x) \in \omega(u, v) \) and for all \( (y, z) \) in the intersection
\[
\left\{ (y, z) \in \mathbb{R}^n \times \mathbb{R}^n: \dist\left( (y, z), \omega(u, v) \right) < \varepsilon \right\} \cap \left\{ (y, z) \in \mathbb{R}^n \times \mathbb{R}^n: H(\bar x, \bar x) < H(y, z) < H(\bar x, \bar x) + \eta \right\}
\]
the following holds:
\[
\varphi'\left( H(y, z) - H(\bar x, \bar x) \right) \| \nabla H(y, z) \| \ \geq \ 1.
\]
Then, according to Theorem \ref{thm}\((viii)\), \( \lim_{t \to +\infty} \dist\left( \big( u(t), v(t) \big), \omega(u, v) \right) = 0 \), therefore, there exists \( t_2 \geq 0 \) such that
\[
\lim_{t \to +\infty} \dist\left( \big( u(t), v(t) \big), \omega(u, v) \right) \ < \ \varepsilon, \text{ for all } t \geq t_2.
\]
Now, by \( \lim_{t \to +\infty} H(u(t), v(t)) = H(\bar x, \bar x) \) and \( H(u(t), v(t)) > H(\bar x, \bar x) \) there exist \( t_3 \geq t_2 \) such that
\[
H(\bar x, \bar x) \ < \ H(u(t), v(t)) \ < \ H(\bar x, \bar x) + \eta, \text{ for all } t \geq t_3.
\]
Hence, we have
\[
\varphi'\left( H(u(t), v(t)) - H(\bar x, \bar x) \right) \| \nabla H(u(t), v(t)) \| \ \geq \ 1 \text{ for all } t \geq t_3.
\]
By Theorem \ref{thm}\((ii)\) and \((v)\) we deduce
\begin{equation}\label{eq}
\begin{split}
    \frac{d}{dt} \varphi\left( H(u(t), v(t)) - H(\bar x, \bar x) \right) \ &= \ \varphi'\left( H(u(t), v(t)) - H(\bar x, \bar x) \right) \frac{d}{dt} H(u(t), v(t)) \\
    &\leq \ \frac{ -A(t) \|\dot{x}(t)\|^2 - B(t) \|\ddot{x}(t)\|^2 }{q(t) \| \ddot x(t) \| + p(t) \| \dot x(t) \|}
\end{split}
\end{equation}
for all \( t \geq t_3 \). By integrating the last inequality on \( [t_3, t] \), we deduce
\[
\varphi\left( H(u(t), v(t)) \right) - H(\bar x, \bar x) + \int_{t_3}^t \frac{A (s) \| \dot x(s) \|^2 + B (s) \|\ddot{x}(s)\|^2}{q(s) \| \ddot x(s) \| + p(s) \| \dot x(s) \|} ds \ \leq \ \varphi\left( H(u(t_3), v(t_3)) \right) - H(\bar x, \bar x).
\]
Since for all \( s \geq t_3 \) \( \varphi, \a, \b, p, q, A \) and \( B \) are non-negative,  and \( t \) was arbitrary chosen, we deduce
\[
0 \ \leq \ \int_{t_3}^{+\infty} \frac{B(s) \| \ddot x(s) \|^2 + A (s) \| \dot x(s) \|^2}{q(s) \| \ddot x(s) \| + p(s) \| \dot x(s) \|} ds \ \leq \ \varphi\left( H(u(t_3), v(t_3)) \right) - H(\bar x, \bar x),
\]
which leads to
\[
\frac{A(t) \| \dot x(t) \|^2+B(t) \| \ddot x(t) \|^2}{q(t) \| \ddot x(t) \| + p(t) \| \dot x(t) \|} \in L^1 \big( [t_0, +\infty), \mathbb{R} \big).
\]

Now according to Lemma \ref{fosikacska} there exists $a,b>0$ such that for all sufficiently large $t$ one has
$$\frac{A(t) \| \dot x(t) \|^2+B(t) \| \ddot x(t) \|^2}{q(t) \| \ddot x(t) \| + p(t) \| \dot x(t) \|} \ge a\| \dot x(t) \|+b\g(t)\| \ddot x(t) \|,$$
hence
\[
\| \dot x(t) \|, \gamma(t) \| \ddot x(t) \| \in L^1 \big( [t_0, +\infty), \mathbb{R} \big).
\]
Now, since  $\dot{x}(t)\in L^1([t_0,+\infty),\R^n)$, it follows that there exists and is finite the limit $\lim_{t\To+\infty}x(t).$ But, obviously, then one must have
$$\lim_{t\To+\infty}x(t)=\ol x.$$
\end{proof}

\begin{corollary}\label{sa}Assume that $g$ is semi-algebraic and bounded from below and for $u_0,v_0\in\R^n$, let $x$  be the unique strong global solution of \eqref{dysy}. If $x$ is bounded then all the conclusions of Theorem \ref{convergence} hold.
\end{corollary}
\begin{proof}
  We need only to prove that $H$ is a KL function. Since the class of semi-algebraic functions is closed under addition (see for example \cite{b-sab-teb}) and
$(x,y) \mapsto \frac12\|x-y\|^2$ is semi-algebraic, we obtain that the  function
$$H:\R^m\times\R^m\To\R\cup\{+\infty\},\,H(x,y)=g(x)+\frac12\|y-x\|^2$$ is semi-algebraic. Consequently, $H$ is a KL function.
\end{proof}

\begin{remark}\label{cond-x-bound}
Note that in the hypotheses of Theorem \ref{convergence} one needs to assume that the  trajectory $x$ generated by \eqref{dysy} is bounded and $g$ is  bounded from below.

In order to validate these assumptions it is enough to assume that  $g$ is coercive,
that is $$\lim_{\|u\|\rightarrow+\infty}g(u)=+\infty.$$

Indeed, notice that $g$ is bounded from below, being a continuous and coercive function (see for example \cite{rock-wets}). From \eqref{ineqfos} it follows that $\g(t)\dot{x}(t)+x(t)$ is contained for every $t\geq t_1$ in a lower level set
of $g$, which is a bounded set due to the coercivity assumption. But, according to Remark \ref{abstgamma}, $\g(t)\dot{x}(t)\To 0,\,t\To+\infty$, hence one can easily deduce that $x$  is bounded.
\end{remark}

\section{Convergence rates}
We show convergence rates of the generated trajectories and function values in terms of the {\L}ojasiewicz exponent of $H.$
Let  $C(t)=\max\left\{\frac{q(t)}{\gamma(t)},\,p(t)\right\}$ and  {for the needs of this section} assume that
\[{\bf (L10)}\,\,\,\,\int_{t_0}^{+\infty}\frac{1}{C(t)}dt=+\infty.\]

The next theorem establishes the convergence rates {of the generated trajectories}.

\begin{theorem}\label{thmrates}

    Assume that the conditions $({\bf L1})-({\bf L10})$ hold. Let \(x\) be a unique strong global solution of \eqref{dysy}. Consider the following regularization of \(g\):
    \[
    H: \mathbb{R}^n \times \mathbb{R}^n \to \mathbb{R}, \ H(x, y) = g(x) + \frac{1}{2} \| x - y \|^2
    \]
    and suppose that \(x\) is bounded. Let \( \bar x \in \crit(g) \) be such that \( \lim_{t \to +\infty} x(t) = \bar x \) and suppose that \(H\) fulfills the \L ojasiewicz property at \( (\bar x, \bar x) \) with the exponent \(\theta\) {(see Definition \ref{Lexpo} in the Appendix B}). Then there exist \( a,R, t_4 > 0 \) such that for all \( t \geq t_4 \) following holds:

    \begin{itemize}

        \item[(i)] If \( \theta \in \left( 0, \frac{1}{2} \right) \), then \(x\) converges in finite time; Further, if $\bar x$ is a local minimizer then $g(x(t)) = g(\bar{x})$ in finite time;

        \item[(ii)] If \( \theta = \frac{1}{2} \), then \( \| x(t) - \bar x \| \ \leq \ a e^{-\int_{t_4}^t \frac{R}{C(s)} ds} \); Further,  $g(x(t)) - g(\bar{x})\le \frac{L_g a^2}{2} e^{-2\int_{t_4}^t \frac{R}{C(s)}\,ds}$

        \item[(iii)] If \( \theta \in \left( \frac{1}{2}, 1 \right) \), then \( \| x(t) - \bar x \| \ \leq \ a\left( \int_{t_4}^t \frac{ds}{C(s)} \right)^{-\frac{1 - \theta}{2 \theta - 1}} \); Further, $g(x(t)) - g(\bar{x})\le \frac{L_g a^2}{2} \left(\int_{t_4}^t \frac{ds}{C(s)}\right)^{-\frac{2(1-\theta)}{2\theta - 1}}$.

    \end{itemize}

\end{theorem}

\begin{proof}

Let $t_2\ge t_1$ such that all the conditions $({\bf L1})-({\bf L10})$ hold. Define for every $t \geq t_2$ the following functional
\[
\sigma(t) \ = \ \int_{t}^{+\infty} \left( \| \dot x(s) \| + \gamma(s) \| \ddot x(s) \| \right) ds.
\]
For a fixed $t \in [t_2, +\infty)$ and $T \geq t$ we have
\[
\| x(t) - \bar x \| \ = \ \left\| x(T) - \bar x - \int_t^T \dot x(s) ds \right\| \ \leq \ \| x(T) - \bar x \| + \int_t^T \| \dot x(s) \| ds.
\]
Taking the limit as $T \to +\infty$ we deduce
\begin{equation}\label{xsigma}
    \| x(t) - \bar x \| \ \leq \ \int_t^{+\infty} \| \dot x(s) \| ds \ \leq \ \sigma(t).
\end{equation}
Notice that according to Lemma \ref{fosikacska} one has for all $t\ge t_2$ that

\begin{equation}\label{eq_1}
    \frac{B(t) \| \ddot x(t) \|^2 + A(t) \| \dot x(t) \|^2}{q(t) \| \ddot x(t) \| + p(t) \| \dot x(t) \|} \ \geq \ m \Big( \| \dot x(t) \| + \gamma(t) \| \ddot x(t) \| \Big)
\end{equation}
where $m \ = \ \min(a,b) > \ 0.$

As we have shown before in the proof of Theorem \ref{convergence}, if there exists \( t_1' \geq t_1 \) such that \( H(u(\bar t), v(\bar t)) = H(\bar x, \bar x) \), then \(x\) is constant on \( [t_1', +\infty) \), and the conclusion of the theorem follows.

Now, if \( H(u(t), v(t)) > H(\bar x, \bar x) \) for every \( t \geq t_1 \), then according to Theorem \ref{thm}\((i)\) and \( \lim_{t \to +\infty} x(t) = \bar x \), there exist \( \varepsilon > 0 \) and \( t_3 \geq t_2 \) such that \( \| (u(t), v(t)) - (\bar x, \bar x) \| < \varepsilon \) for all \( t \geq t_3 \). Since \(H\) fulfills the \L ojasiewicz property at \( (\bar x, \bar x) \in \crit H \) with the exponent \(\theta\), by \eqref{eq} we have
\[
K \frac{d}{dt} \left( H(u(t), v(t)) - H(\bar x, \bar x) \right)^{1 - \theta} \ \leq \ \frac{ -A(t) \|\dot{x}(t)\|^2 - B(t) \|\ddot{x}(t)\|^2 }{q(t) \| \ddot x(t) \| + p(t) \| \dot x(t) \|}, \text{ for all } t \geq t_3.
\]
By \eqref{eq_1} we deduce for \( M = \frac{m}{K} \)
\[
M \left( \| \dot x(t) \| + \gamma(t) \| \ddot x(t) \| \right) + \frac{d}{dt} \left( H(u(t), v(t)) - H(\bar x, \bar x) \right)^{1 - \theta} \ \leq \ 0, \text{ for all } t \geq t_3.
\]
After integrating the last relation on \( [t, t_4] \), for \( t_4> t > t_3 \), we obtain
\[
M \int_t^{t_4} \left( \| \dot x(s) \| + \gamma(s) \| \ddot x(s) \| \right) ds + \left( H(u(t_4), v(t_4)) - H(\bar x, \bar x) \right)^{1 - \theta} \ \leq \ \left( H(u(t), v(t)) - H(\bar x, \bar x) \right)^{1 - \theta}.
\]
Letting \( t_4 \to +\infty \) we deduce
\[
M \sigma(t) \ \leq \ \left( H(u(t), v(t)) - H(\bar x, \bar x) \right)^{1 - \theta}, \text{ for all } t \geq t_4.
\]
On the other hand,
\[
\left( H(u(t), v(t)) - H(\bar x, \bar x) \right)^\theta \ \leq \ K \| \nabla H(u(t), v(t)) \|
\]
and by Theorem \ref{thm}\((v)\)
\[
\| \nabla H(u(t), v(t)) \| \ \leq \ q(t) \| \ddot x(t) \| + p(t) \| \dot x(t) \|, \text{ for all } t \geq t_4.
\]
Thus,
\[
M \sigma(t) \ \leq \ K^\frac{1 - \theta}{\theta} \Big( q(t) \| \ddot x(t) \| + p(t) \| \dot x(t) \| \Big)^\frac{1 - \theta}{\theta}, \text{ for all } t \geq t_4.
\]
Notice that deploying $C(t)=\max\left\{\frac{q(t)}{\gamma(t)},\,p(t)\right\}$ we get
\[
q(t) \| \ddot x(t) \| + p(t) \| \dot x(t) \| \ \leq C(t) \Big( \gamma(t) \| \ddot x(t) \| + \| \dot x(t) \| \Big), \text{ for all } t \geq t_4,
\]
and thus
\[
M \sigma(t) \ \leq \ \left(  K C(t) \right)^\frac{1 - \theta}{\theta} \Big( \gamma(t) \| \ddot x(t) \| + \| \dot x(t) \| \Big)^\frac{1 - \theta}{\theta}, \text{ for all } t \geq t_4.
\]
Finally, using \( -\dot \sigma(t) = \gamma(t) \| \ddot x(t) \| + \| \dot x(t) \| \) we obtain for \( R = \frac{M^\frac{\theta}{1 - \theta}}{K} > 0 \)
\begin{equation}\label{key_ineq}
    -\frac{R}{C(t)} \sigma^\frac{\theta}{1 - \theta}(t) \ \geq \ \dot \sigma(t), \text{ for all } t \geq t_4,
\end{equation}
or equivalently   for all $t \geq t_4,$ holds
\begin{equation}\label{key_ineq1}
    -\frac{R}{C(t)} \geq \ \dot \sigma(t)\sigma^\frac{\theta}{\theta-1}(t)=\left\{\begin{array}{lll} \frac{d}{dt}\ln(\sigma(t),\text{ if }\theta=\frac12\\
    \\
    \frac{\theta-1}{2\theta-1}\frac{d}{dt}\sigma^\frac{2\theta-1}{\theta-1}(t),\text{ if }\theta\neq\frac12.\end{array}\right.
\end{equation}

\textbf{Case \( \theta = \frac{1}{2} \).}

In this case \eqref{key_ineq1} becomes
\[
\frac{d}{dt}\ln(\sigma(t) \ \leq \  -\frac{R}{C(t)} \ \forall t \geq t_4.
\]
Integrating on \( [t_4, T] \) we deduce
\[
\ln\frac{\sigma(T)}{\sigma(t_4)}  \ \leq \ -\int_{t_4}^T \frac{R}{C(s)} ds. \]
In other words,  for all $t\ge t_4$ one has \( \sigma(t)\le \sigma(t_4) \left( e^{-\int_{t_4}^t \frac{R}{C(s)} ds} \right)\).
Using \eqref{xsigma} and letting $a=\sigma(t_4)$ we obtain
\[
\| x(t) - \bar x \| \ \leq \ a e^{-\int_{t_4}^t \frac{R}{C(s)} ds} \ \forall t \geq t_4.
\]

\textbf{Case \( 0 \leq \theta < \frac{1}{2} \).} We show by contradiction. Assume that \( \sigma(t) \neq 0 \) for all \( t \geq t_4 \). Then,
\[
\frac{d}{dt} \left( \sigma^\frac{2 \theta-1}{\theta-1}(t) \right) \ \leq \  -\frac{2 \theta-1}{\theta-1} \frac{R}{C(t)} \ \forall t \geq t_4.
\]
By integrating on \( [t_4, t] \) we deduce
\[
\sigma^\frac{1 - 2 \theta}{1 - \theta}(t) \ \leq \ \sigma^\frac{1 - 2 \theta}{1 - \theta}(t_4)-\frac{1 - 2 \theta}{1 - \theta} R \int_{t_4}^t \frac{ds}{C(s)} \ \forall t \geq t_4
\]
and since according to $({\bf L10})$ the right hand side goes to $-\infty$ as $t\to+\infty$ we get that \( \sigma(t) < 0 \) for all sufficiently large $t$, which is a contradiction with \( \sigma(t) > 0 \ \forall t \geq t_4 \). That means that there exists \( T \geq t_4 \) such that \( \sigma(T) = 0 \). Since \( \sigma \) is nonnegative and nonincreasing, \( \sigma(t) = 0 \) for all \( t \geq T \). Therefore,
\eqref{xsigma} implies that \( x \) is constant on \( [T, +\infty) \), i.e. $x(t)=\bar x$ for all \( t\in[T, +\infty) \).

\textbf{Case \( \frac{1}{2} < \theta < 1 \).}
Using \eqref{key_ineq1} we deduce
\[
\frac{d}{dt} \left( \sigma^{-\frac{2 \theta-1}{1 - \theta}(t)} \right) \ \geq \  \frac{2 \theta - 1}{1 - \theta} \frac{R}{C(t)} \ \forall t \geq t_4.
\]
Again, by integrating on \( [t_4, t] \) we obtain
\[
\sigma^{-\frac{2 \theta-1}{1 - \theta}(t)} \ \geq \sigma^{-\frac{2 \theta-1}{1 - \theta}}(t_4)+ \ \frac{2 \theta - 1}{1 - \theta} R \int_{t_4}^t \frac{ds}{C(s)} \ \forall t \geq t_4
\]
or
\[
\sigma(t) \ \leq \ \left(\sigma^{-\frac{2 \theta-1}{1 - \theta}}(t_4)+ \frac{2 \theta - 1}{1 - \theta} R \int_{t_4}^t \frac{ds}{C(s)} \right)^{-\frac{1 - \theta}{2 \theta - 1}} \ \forall t \geq t_4.
\]
Using \eqref{xsigma} we finally deduce
\[
\| x(t) - \bar x \| \ \leq \ \left( \frac{2 \theta - 1}{1 - \theta} R \int_{t_4}^t \frac{ds}{C(s)} \right)^{-\frac{1 - \theta}{2 \theta - 1}} \ \forall t \geq t_4.
\]

In order to obtain the corresponding rates for the function values it is enough to deploy the descent lemma.  Indeed, since $\nabla g$ is $L_g$-Lipschitz continuous and $\nabla g(\bar{x}) = 0$, the descent lemma yields
\[
g(x(t)) - g(\bar{x}) \le \langle \nabla g(\bar{x}), x(t) - \bar{x} \rangle + \frac{L_g}{2}\|x(t) - \bar{x}\|^2 = \frac{L_g}{2}\|x(t) - \bar{x}\|^2.
\]
 For (i) note that $g(x(t)) \ge g(\bar{x})$ holds near $\bar{x}$ whenever $\bar{x}$ is a local minimizer. The conclusions of (ii) and (iii) follow from the rates established for $\|x(t) - \bar{x}\|.$
\end{proof}

As we have seen we need several conditions in order to obtain our results. In what follows we show that these conditions imply a strong corelation among the parameters $\a,\b$ and $\g$.

\begin{remark}\label{Asymptotic compatibility of the parameters}
The assumptions considered in this paper yield the following
asymptotic relations between the coefficients $\alpha,\b$ and $\g$:
\[
\begin{array}{c|c}
\textbf{Assumptions} & \textbf{Consequence}\\
\hline
\textbf{(L5)--(L7)}
& \beta(t)=\Omega\big(\alpha^2(t)\big)\\[2mm]
\textbf{(L1), (L4), (L7)}
& \beta(t)=O\big(\alpha^2(t)\big)\\[2mm]
\textbf{(L1), (L4)--(L7)}
& \beta(t)=\Theta\big(\alpha^2(t)\big)
\end{array}
\]
Here, all asymptotic relations are understood as
$t\to+\infty$.

We emphasize that the last conclusion does not require the
entire condition \textbf{(L1)}: its second inequality,
namely
\(
\frac{\dot{\beta}(t)}{\beta(t)}
\leq K_2\alpha(t),
\)
is sufficient to establish the upper estimate
$\beta(t)=O(\alpha^2(t))$, together with
\textbf{(L4)} and \textbf{(L7)}.
Further, \textbf{(L7)} implies 
\[
\gamma(t)=\Theta\left(\frac{1}{\alpha(t)}\right).
\]
\end{remark}
Accordingly to the corelations among $\a,\b,\g$ emphasized in the previous remark we obtain superlinear rates as a corollary of Theorem \ref{thmrates}.
\begin{corollary}[Superlinear convergence rates]\label{cor:superlinear}
Assume that $g$ is coercive and
$H : \mathbb{R}^n \times \mathbb{R}^n \to \mathbb{R}$,
$H(x,y)=g(x)+\frac12\|x-y\|^2$, is a KL function. Fix $\delta>0$ and $0<\kappa<\frac{\sqrt{1+2\delta}-1}{L_g\delta^2}$. (Note that $\delta=1$, this reduces to
$0<\kappa<\frac{\sqrt{3}-1}{L_g}$.)
Let $\alpha:[t_0,+\infty)\to(0,+\infty)$ be a nondecreasing
$C^2$ function such that
\begin{equation}\label{eq:alpha_cond}
\lim_{t\to+\infty}\alpha(t)=+\infty,\qquad
\lim_{t\to+\infty}\frac{\dot{\alpha}(t)}{\alpha^2(t)}=0,
\end{equation}
and
\begin{equation}\label{eq:alpha_derivative_cond}
\lim_{t\to+\infty}\alpha(t)\frac{d}{dt}
\left(\frac{\delta+1-\delta\dot{\alpha}(t)/\alpha^2(t)}
{2\kappa\alpha^2(t)}\right)=0.
\end{equation}
Set $\beta(t)=\kappa\alpha^2(t)$ and
$\gamma(t)=\delta/\alpha(t)$.

Let $x$ be the unique strong global solution of \textbf{(DS)}, which converges to some $\bar{x} \in \mathrm{crit}(g)$, and suppose that $H$ fulfills the \L{}ojasiewicz property at $(\bar{x}, \bar{x})$ with exponent $\theta \in (0, 1)$. Then there exist constants $a, c, t_4 > 0$ such that for all $t \ge t_4$ the following hold:
\begin{itemize}
    \item[\textbf{(i)}] If $\theta \in \left(0, \frac{1}{2}\right)$, then $x$ converges to $\bar{x}$ in finite time;
    \item[\textbf{(ii)}] If $\theta = \frac{1}{2}$, then $x(t)$ converges at the superlinear-exponent rate
    \[
    \|x(t) - \bar{x}\| \le a e^{-c \int_{t_4}^t \alpha(s)\,ds};
    \]
    \item[\textbf{(iii)}] If $\theta \in \left(\frac{1}{2}, 1\right)$, then
    \[
    \|x(t) - \bar{x}\| \le a \left(\int_{t_4}^t \alpha(s)\,ds\right)^{-\frac{1-\theta}{2\theta - 1}}.
    \]
\end{itemize}
\end{corollary}

\begin{proof}
By the definition of $\beta$ and $\gamma$, $\beta$ is continuously differentiable and non-decreasing with $\frac{\dot{\beta}(t)}{\beta(t)} = \frac{2\dot{\alpha}(t)}{\alpha(t)}$, while $\dot{\gamma}(t) = -\frac{\delta \dot{\alpha}(t)}{\alpha^2(t)}$. By \eqref{eq:alpha_cond}, for any $\varepsilon \in (0, 1)$ there exists $T_\varepsilon \ge t_0$ such that $0 \le \frac{\dot{\alpha}(t)}{\alpha^2(t)} \le \varepsilon$ for all $t \ge T_\varepsilon$, which implies $\dot{\gamma}(t) \in [-\delta\varepsilon, 0]$ and $\lim_{t \to +\infty} \dot{\gamma}(t) = 0$.

We verify that conditions \textbf{(L1)}--\textbf{(L10)} hold for $t$ sufficiently large:
\begin{itemize}
  \item[  ] \textbf{(L1)} We have $\frac{\dot{\alpha}(t) + \alpha(t)\dot{\beta}(t)/\beta(t)}{\beta(t)} = \frac{3\dot{\alpha}(t)}{\kappa \alpha^2(t)} \le \frac{3\varepsilon}{\kappa} \le K_1 \kappa \alpha^2(t) = K_1\beta(t)$ and $\frac{\dot{\beta}(t)}{\beta(t)} = \frac{2\dot{\alpha}(t)}{\alpha(t)} \le 2\varepsilon \alpha(t) \le K_2 \alpha(t)$ for $t$ large enough.
 \item[ ]  \textbf{(L3)-(L4)}Since $\alpha(t)\gamma(t) + \dot{\gamma}(t) + 1 = \delta + 1 + \dot{\gamma}(t) \to \delta + 1 > 0$, \textbf{(L3)} holds for large $t$. Moreover, $\frac{1}{\gamma(t)} = \frac{1}{\delta}\alpha(t)$, so \textbf{(L4)} holds with $C_1 = \frac{1}{\delta}$.
    \item[] \textbf{(L5)} Direct substitution yields
    \[
    \frac{2\dot{\alpha}(t)\alpha(t)\beta(t) - 2\alpha^2(t)\dot{\beta}(t)}{\beta^2(t)} + \frac{\alpha^4(t)}{\beta^2(t)\gamma(t)} = -\frac{2\dot{\alpha}(t)}{\kappa \alpha(t)} + \frac{1}{\kappa^2 \delta}\alpha(t) \le \frac{1}{\kappa^2 \delta}\alpha(t),
    \]
    so \textbf{(L5)} holds with $C_2 = \frac{1}{\kappa^2\delta}$.
    \item[]\textbf{(L6)} Similarly,
    \[
    \frac{\big(\dot{\alpha}(t)\beta(t) + \alpha(t)\dot{\beta}(t) + L_g \beta^2(t)\big)^2}{\gamma(t)\beta^4(t)} = \frac{\alpha(t)}{\delta\kappa^2}\left(\frac{3\dot{\alpha}(t)}{\alpha^2(t)} + L_g\kappa\right)^2 \le \frac{(3\varepsilon + L_g\kappa)^2}{\delta\kappa^2}\alpha(t),
    \]
    confirming \textbf{(L6)} with $C_3 = \frac{(3\varepsilon + L_g\kappa)^2}{\delta\kappa^2}$.
    \item[] \textbf{(L7)--(L8)} Since $\frac{\alpha(t)}{\beta(t)} = \frac{1}{\kappa \alpha(t)} = \frac{1}{\kappa\delta}\gamma(t)$, \textbf{(L7)} holds with $C_4 = C_5 = \frac{1}{\kappa\delta}$. Furthermore, for $t \ge T_\varepsilon$ we have $\dot{\gamma}(t) \le 0$, hence
    \[
    \sqrt{\frac{\alpha(t)\gamma(t) + \dot{\gamma}(t) + 1}{\beta(t)}} \le \frac{\sqrt{\delta + 1}}{\sqrt{\kappa}\,\alpha(t)} = \sqrt{\kappa(\delta+1)}\,\frac{\alpha(t)}{\beta(t)},
    \]
    which gives \textbf{(L8)} with $C_6 = \sqrt{\kappa(\delta+1)}$.
    \item[] \textbf{(L2) \& (L9)} Notice that \textbf{(L9)} implies \textbf{(L2)}. Multiplying the right-hand side of \textbf{(L9)} by $\alpha(t)$ and letting $t \to +\infty$ (noting that $\dot{\gamma}(t) \to 0$ and $\big(\frac{\alpha(t)\gamma(t)+\dot{\gamma}(t)+1}{2\beta(t)}\big)' \to 0$), we obtain the asymptotic limit
    \[
    \frac{1 - C_7}{\kappa} - \frac{L_g^2\delta^3\kappa}{2} - L_g\delta(1 + C_7) - 2C_7\sqrt{\frac{\delta+1}{\kappa}}.
    \]
    Meanwhile, multiplying the left-hand side of \textbf{(L9)} by $\alpha(t)$ and taking $t \to +\infty$ yields
    \[
    \frac{C_7}{2 - 4C_8}\frac{\kappa}{\delta}\left(C_8\delta\left(L_g\delta + \frac{1}{\kappa} + 2\sqrt{\frac{\delta+1}{\kappa}}\right)\right)^2.
    \]
    Since $0 < \kappa < \frac{\sqrt{1+2\delta}-1}{L_g\delta^2}$ is equivalent to $\frac{1}{\kappa} - \frac{L_g^2\delta^3\kappa}{2} - L_g\delta > 0$, by choosing $C_7 > 0$ and $C_8 \in \left(0, \frac{1}{2}\right)$ sufficiently small, the asymptotic limit of the right-hand side strictly exceeds that of the left-hand side (and is strictly positive). Hence \textbf{(L9)}, and consequently \textbf{(L2)}, hold for all $t$ sufficiently large.
    \item[] \textbf{(L10)} By definition of $q(t)$ and $p(t)$ in Theorem~\ref{thm}\textbf{(v)}, for large $t$ we have
    \[
    \frac{q(t)}{\gamma(t)} = \frac{1}{\beta(t)\gamma(t)} = \frac{1}{\kappa\delta \alpha(t)}
    \]
    and
    \[
    p(t) = L_g\gamma(t) + \frac{\alpha(t)}{\beta(t)} + 2\sqrt{\frac{\alpha(t)\gamma(t)+\dot{\gamma}(t)+1}{\beta(t)}} \le \left(L_g\delta + \frac{1}{\kappa} + 2\sqrt{\frac{\delta+1}{\kappa}}\right)\frac{1}{\alpha(t)}.
    \]
    Thus, there exists $M_0 > 0$ such that $C(t) = \max\left\{\frac{q(t)}{\gamma(t)},\,p(t)\right\} \le \frac{M_0}{\alpha(t)}$ for all $t$ large enough. Since $\alpha(t) \to +\infty$, we have $\int_{t_0}^{+\infty} \frac{1}{C(s)}\,ds = +\infty$, so \textbf{(L10)} holds.
\end{itemize}
Finally, by coercivity of $g$, Remark~\ref{cond-x-bound} ensures that $g$ is bounded from below and the trajectory $x$ is bounded, so Theorem~\ref{convergence} guarantees $\lim_{t \to +\infty} x(t) = \bar{x} \in \mathrm{crit}(g)$. Applying Theorem~\ref{thmrates} together with the lower bound $\frac{1}{C(s)} \ge \frac{1}{M_0}\alpha(s)$ immediately yields conclusions \textbf{(i)}--\textbf{(iii)} with $c = \frac{R}{M_0}$.
\end{proof}

\begin{remark}\label{ex:rates}
Corollary~\ref{cor:superlinear} covers several explicit families of unbounded damping coefficients:
\begin{enumerate}
    \item \textbf{Polynomial damping (compressed/Gaussian super-exponential rates):} Let $\alpha(t) = \alpha_0 t^p$ and $\beta(t) = \beta_0 t^{2p}$ for $t \ge t_0 > 0$, where $p > 0$ and $0 < \beta_0 < \frac{\sqrt{3}-1}{L_g}\alpha_0^2$. Then $\frac{\dot{\alpha}(t)}{\alpha^2(t)} = \frac{p}{\alpha_0}t^{-p-1} \to 0$ as $t \to +\infty$. For $\theta = \frac{1}{2}$, we obtain the superlinear-exponent convergence rate
    \[
    \|x(t) - \bar{x}\| \le a e^{-\tilde{c} t^{p+1}},
    \]
    which reduces to the Gaussian rate $\mathcal{O}\big(e^{-\tilde{c} t^2}\big)$ for $p = 1$ , and achieves arbitrarily high polynomial powers $t^{p+1}$ in the exponent for $p > 1$. For $\theta \in \left(\frac{1}{2}, 1\right)$, we obtain $\|x(t) - \bar{x}\| = \mathcal{O}\big(t^{-(p+1)\frac{1-\theta}{2\theta-1}}\big)$.
    
    \item \textbf{Exponential damping (double-exponential rates):} Let $\alpha(t) = \alpha_0 e^{\lambda t}$ and $\beta(t) = \beta_0 e^{2\lambda t}$ for $t \ge t_0 \ge 0$, where $\lambda > 0$ and $0 < \beta_0 < \frac{\sqrt{3}-1}{L_g}\alpha_0^2$. Then $\frac{\dot{\alpha}(t)}{\alpha^2(t)} = \frac{\lambda}{\alpha_0}e^{-\lambda t} \to 0$ as $t \to +\infty$. Consequently, for $\theta = \frac{1}{2}$ the trajectory converges at the double-exponential rate
    \[
    \|x(t) - \bar{x}\| \le a e^{-\tilde{c} e^{\lambda t}},
    \]
    while for $\theta \in \left(\frac{1}{2}, 1\right)$ we still obtain exponential convergence $\|x(t) - \bar{x}\| = \mathcal{O}\big(e^{-\lambda\frac{1-\theta}{2\theta-1}t}\big)$.

    \item \textbf{Logarithmic damping:} Let
    $\alpha(t)=\alpha_0\log(e+t)$,
    $\beta(t)=\kappa\alpha^2(t)$, and
    $\gamma(t)=\delta/\alpha(t)$, where $\alpha_0,\delta>0$ and
    $0<\kappa<(\sqrt{1+2\delta}-1)/(L_g\delta^2)$.
    Then $\alpha(t)\to+\infty$,
    $\dot{\alpha}(t)/\alpha^2(t)\to0$, and
    \eqref{eq:alpha_derivative_cond} holds. For $\theta=1/2$,
    \[
    \|x(t)-\bar{x}\|\leq
    a\exp\!\left(-\tilde c\int_{t_4}^t\log(e+s)\,ds\right),
    \]
    which has an exponent of order $t\log t$.
    
    \item \textbf{Mixed polynomial--logarithmic--exponential damping:}
    The same conditions hold for
    \[
    \alpha(t)=\alpha_0(1+t)^p[\log(e+t)]^r e^{\lambda t},
    \quad \beta(t)=\kappa\alpha^2(t),\quad
    \gamma(t)=\frac{\delta}{\alpha(t)},
    \]
    where $p,r,\lambda\geq0$, at least one of $p,r,\lambda$ is
    positive, and the constants satisfy the preceding restriction.
    In these cases \eqref{eq:alpha_derivative_cond} is satisfied.
    Since the proof of the corollary gives
    $C(t)\leq M_0/\alpha(t)$ for all sufficiently large $t$,
    condition \textbf{(L10)} follows from
    $\int_{t_4}^\infty C(t)^{-1}\,dt\geq
    M_0^{-1}\int_{t_4}^\infty\alpha(t)\,dt=+\infty$.
\end{enumerate}

All the parameter choices considered above satisfy the
same asymptotic scaling between the coefficients:
\[
\beta(t)=\Theta(\alpha^2(t)),
\qquad
\gamma(t)=\Theta\left(\frac{1}{\alpha(t)}\right),
\qquad t\to+\infty.
\]
Indeed, in all four cases there exist positive constants
$\kappa,\delta$ such that
\[
\beta(t)=\kappa\alpha^2(t),
\qquad
\gamma(t)=\frac{\delta}{\alpha(t)}.
\]
Moreover, each choice satisfies
\[
\lim_{t\to+\infty}\frac{\dot\alpha(t)}{\alpha^2(t)}=0,
\qquad
\int_{t_0}^{+\infty}\alpha(t)\,dt=+\infty.
\]

Consequently, in the critical case $\theta=\frac12$,
the convergence rate is governed by
\[
\|x(t)-\bar x\|
\leq a\exp\left(-c\int_{t_4}^t\alpha(s)\,ds\right).
\]
Thus, polynomial, logarithmic, exponential and mixed
damping profiles lead to different superlinear-exponent
convergence rates, despite sharing the same asymptotic
scaling between $\alpha$, $\beta$ and $\gamma$.
\end{remark}

\section{On the KL-property of $g$ and the strongly convex case}

According to Theorem 3.6 \cite{LiP}, if $g$ has the KL property with {\L}ojasiewicz exponent $\t\in\left[\frac12,1\right)$ at $\ol x\in\R^m,$ then the function $H:\R^m\times\R^m\To\R\cup\{+\infty\},\,H(x,y)=g(x)+\frac12\|y-x\|^2$ has the KL property at $(\ol x,\ol x)\in\R^m\times\R^m$ with the same  exponent $\t.$ Moreover,  if $g$ has the KL property with {\L}ojasiewicz exponent $\t\in\left(0,\frac12\right)$ at $\ol x\in\R^m,$ then the function $H:\R^m\times\R^m\To\R\cup\{+\infty\},\,H(x,y)=g(x)+\frac12\|y-x\|^2$ has the KL property at $(\ol x,\ol x)\in\R^m\times\R^m$ with the exponent $\t=\frac12.$  For the sake of completness let us show this claim,

\begin{proposition}\label{Exact KL exponent of the regularized function}
\label{prop:exact-KL-exponent}
Let $g:\mathbb{R}^n\to\mathbb{R}$ be continuously
differentiable and suppose that $g$ satisfies the
Kurdyka--\L{}ojasiewicz inequality at a critical point
$\bar{x}\in\mathbb{R}^n$ with an exponent
$\theta\in(0,\frac12)$.

Consider the function
\(
H:\mathbb{R}^n\times\mathbb{R}^n\to\mathbb{R},
\qquad
H(x,y)=g(x)+\frac12\|x-y\|^2.
\)
Then $H$ satisfies the KL inequality at
$(\bar{x},\bar{x})$ with exponent $\frac12$.
Moreover, it cannot satisfy the KL inequality at
$(\bar{x},\bar{x})$ with any exponent smaller than
$\frac12$, and consequently, the exact KL exponent of $H$ at
$(\bar{x},\bar{x})$ is $\frac12$.
\end{proposition}

\begin{proof}
Since $g$ satisfies the KL inequality at $\bar{x}$
with exponent $\theta\in(0,\frac12)$, there exist $c_g,\varepsilon,\eta>0$ such that \(\|\nabla g(x)\|\geq c_g\big(g(x)-g(\bar{x})\big)^\theta
\) whenever \(\|x-\bar{x}\|<\varepsilon, \qquad 0<g(x)-g(\bar{x})<\eta.\)

Without loss of generality, we may assume that $\eta\leq1$. Consequently,
\(
\big(g(x)-g(\bar{x})\big)^\theta
\geq
\big(g(x)-g(\bar{x})\big)^{\frac12},
\)
and hence
\[
\|\nabla g(x)\|
\geq c_g\big(g(x)-g(\bar{x})\big)^{\frac12}.
\]
In other words, $g$ also satisfies the KL inequality at $\bar{x}$ with exponent $\frac12$. According to Theorem 3.6 \cite{LiP} this shows that  $H$ satisfies the KL inequality at
$(\bar{x},\bar{x})$ with exponent $\frac12$.

It remains to show that no smaller exponent is possible.
Suppose, to the contrary, that $H$ satisfies the KL
inequality at $(\bar{x},\bar{x})$ with some exponent
$\widehat{\theta}\in[0,\frac12)$.
Then there exist $c_H,\varepsilon_H,\eta_H>0$
such that
\[
\|\nabla H(x,y)\|
\geq c_H
\big(H(x,y)-H(\bar{x},\bar{x})\big)^{\widehat{\theta}}\mbox{ whenever }
\|(x,y)-(\bar{x},\bar{x})\|<\varepsilon_H,
\qquad
0<H(x,y)-H(\bar{x},\bar{x})<\eta_H.
\]

Consider now the particular choice
\(
x=\bar{x},\qquad y=\bar{x}+h,
\)
where $h\in\mathbb{R}^n\setminus\{0\}$ is sufficiently
small. Since $\nabla g(\bar{x})=0$, we have
\[
H(\bar{x},\bar{x}+h)-H(\bar{x},\bar{x})
=\frac12\|h\|^2\mbox{ and }
\|\nabla H(\bar{x},\bar{x}+h)\|
=\sqrt{2}\,\|h\|.
\]
Consequently,
\(
\sqrt{2}\,\|h\|
\geq
c_H\left(\frac12\|h\|^2\right)^{\widehat{\theta}},
\)
or equivalently,
\(
\sqrt{2}\,\|h\|^{1-2\widehat{\theta}}
\geq c_H2^{-\widehat{\theta}}.
\)
Since $1-2\widehat{\theta}>0$, letting $h\to0$
yields a contradiction.

Thus, $H$ cannot satisfy the KL inequality at
$(\bar{x},\bar{x})$ with any exponent smaller than
$\frac12$. Combining this fact with the first part
of the proof, we conclude that the exact KL exponent
of $H$ at $(\bar{x},\bar{x})$ is $\frac12$.
\end{proof}

 This previous result allows us to reformulate Theorem \ref{thmrates}.

\begin{corollary}\label{cuse}   Assume that $g$ is bounded from below and for $u_0,v_0\in\R^n$, let $x$  be the unique strong global solution of \eqref{dysy}.
Suppose that $x$ is bounded and let $\ol x\in\crit(g)$ be such that $\lim_{t\To+\infty}x(t)=\ol x$ and $g$ fulfills the {\L}ojasiewicz property at $\ol x\in\crit g$ with {\L}ojasiewicz exponent   $\t\in\left(0,1\right)$. If $\t\in\left(0,\frac12\right]$ then the convergence rates stated at {\rm Theorem \ref{thmrates}(ii)} hold,  if  $\t\in\left(\frac12,1\right)$  then the convergence rates stated at  {\rm Theorem \ref{thmrates}(iii)} hold.
\end{corollary}
\begin{proof}
 Indeed, according to Theorem 3.6 \cite{LiP} and Proposition \ref{Exact KL exponent of the regularized function}
$H$ has the {\L}ojasiewicz property at $(\ol x,\ol x)\in\crit H$ with the {\L}ojasiewicz exponent   $\t\in\left[\frac12,1\right)$. Hence, Theorem \ref{thmrates} can be applied.
\end{proof}

In case we assume that the function $g$ is strongly convex, then Theorem \ref{thmrates} assures superlinear convergence rates for the trajectories generated by \eqref{dysy}. We emphasize that the strongly convex function $g$  is coercive, see \cite{BauComb}. According to \cite{rock-wets} the function $g$ is bounded from bellow. Hence, $\omega(x)=\{x^*\}$ and $x(t)\To x^*,\,t\To+\infty.$ According to \cite{attouch-bolte2009}, $g$ satisfies the Kurdyka-{\L}ojasiewicz property at $x^*$ with the {\L}ojasiewicz exponent $\t=\frac12.$  The following result holds.

\begin{corollary}[Strongly convex objective functions]\label{cor:strongly_convex}
Assume that $g : \mathbb{R}^n \to \mathbb{R}$ is $\mu$-strongly convex for some $\mu > 0$ (i.e., $g - \frac{\mu}{2}\|\cdot\|^2$ is convex) with an $L_g$-Lipschitz continuous gradient, and let $\bar{x} \in \mathbb{R}^n$ be the unique global minimizer of $g$. Suppose that the parameters $\alpha, \beta$ satisfy the conditions of Corollary~\ref{cor:superlinear}. Then the regularized function
\[
H : \mathbb{R}^n \times \mathbb{R}^n \to \mathbb{R}, \qquad H(x,y) = g(x) + \frac{1}{2}\|x-y\|^2,
\]
is $\mu_H$-strongly convex with $\mu_H = \frac{\mu + 2 - \sqrt{\mu^2 + 4}}{2} > 0$, and fulfills the \L{}ojasiewicz property at $(\bar{x}, \bar{x})$ with exponent $\theta = \frac{1}{2}$. Consequently, for every initial condition $(u_0, v_0) \in \mathbb{R}^n \times \mathbb{R}^n$, there exist constants $a, c, t_4 > 0$ such that for all $t \ge t_4$:
\[
\|x(t) - \bar{x}\| \le a e^{-c \int_{t_4}^t \alpha(s)\,ds} \qquad \text{and} \qquad 0 \le g(x(t)) - \min g \le \frac{L_g a^2}{2} e^{-2c \int_{t_4}^t \alpha(s)\,ds}.
\]
\end{corollary}

\begin{proof}
Since $g$ is $\mu$-strongly convex, it is coercive and admits a unique global minimizer $\bar{x} \in \mathrm{crit}(g)$ with $g(\bar{x}) = \min_{\mathbb{R}^n} g$. To show that $H$ is strongly convex, observe that for any $(x,y) \in \mathbb{R}^n \times \mathbb{R}^n$,
\[
H(x,y) - \frac{\mu_H}{2}\big(\|x\|^2 + \|y\|^2\big) = \left(g(x) - \frac{\mu}{2}\|x\|^2\right) + \frac{1}{2}\Big((\mu - \mu_H + 1)\|x\|^2 - 2\langle x, y \rangle + (1 - \mu_H)\|y\|^2\Big).
\]
The first term is convex by the $\mu$-strong convexity of $g$. The second term is a quadratic form associated with the $2 \times 2$ block matrix
\[
M = \begin{pmatrix} (\mu - \mu_H + 1)I_n & -I_n \\ -I_n & (1 - \mu_H)I_n \end{pmatrix},
\]
whose eigenvalues are $\lambda_{\pm} = \frac{\mu + 2 \pm \sqrt{\mu^2 + 4}}{2} - \mu_H$. Choosing $\mu_H = \frac{\mu + 2 - \sqrt{\mu^2 + 4}}{2} > 0$ ensures that $\lambda_- = 0$ and $\lambda_+ = \sqrt{\mu^2 + 4} > 0$, so $M$ is positive semidefinite. Thus, $H - \frac{\mu_H}{2}\|\cdot\|^2$ is convex, meaning $H$ is $\mu_H$-strongly convex on $\mathbb{R}^n \times \mathbb{R}^n$.

In particular, $H$ is a KL function  and satisfies the \L{}ojasiewicz inequality globally around its unique minimizer $(\bar{x}, \bar{x})$ with $H(\bar{x}, \bar{x}) = g(\bar{x})$, that is
\[
H(x,y) - H(\bar{x}, \bar{x}) \le \frac{1}{2\mu_H} \|\nabla H(x,y)\|^2 \qquad \forall (x,y) \in \mathbb{R}^n \times \mathbb{R}^n,
\]
which shows that $H$ satisfies Definition \ref{Lexpo} at $(\bar{x}, \bar{x})$ with exponent $\theta = \frac{1}{2}$ and constant $K = \frac{1}{\sqrt{2\mu_H}}$. Applying Corollary~\ref{cor:superlinear}\textbf{(ii)} and Theorem \ref{thmrates}\textbf{(ii)} directly yields the claimed superlinear convergence rates for both $\|x(t) - \bar{x}\|$ and $g(x(t)) - \min g$.
\end{proof}

\appendix

\section{On the local absolute continuity}

In what follows we present the notion that we need in our existence/uniqueness result.

\begin{definition}\label{abs-cont} \rm (see, for instance, \cite{att-sv2011, abbas-att-sv}) A function $x:[0,+\infty)\rightarrow \R^n$  is said to be locally absolutely continuous, if is absolutely continuous on every interval $[0,T],\,T>0$, that is, one of the
following equivalent properties holds:

(i)  there exists an integrable function $y:[0,T]\rightarrow \R^n$ such that $$x(t)=x(0)+\int_0^t y(s)ds \ \ \forall t\in[0,T];$$

(ii) $x$ is continuous and its distributional derivative is Lebesgue integrable on $[0,T]$;

(iii) for every $\varepsilon > 0$, there exists $\eta >0$ such that for any finite family of intervals $I_k=(a_k,b_k)\subseteq [0,T]$ we have the implication
$$\left(I_k\cap I_j=\emptyset \mbox{ and }\sum_k|b_k-a_k| < \eta\right)\Longrightarrow \sum_k\|x(b_k)-x(a_k)\| < \varepsilon.$$
\end{definition}

\begin{remark}\label{rem-abs-cont}\rm (a) It follows from the definition that an absolutely continuous function is differentiable almost
everywhere, its derivative coincides with its distributional derivative almost everywhere and one can recover the function from its derivative $\dot x=y$ by the integration formula (i).

(b) If $x:[0,T]\rightarrow \R^n$ (where $T>0$) is absolutely continuous and $B:\R^n\rightarrow \R^n$ is $L$-Lipschitz continuous
(where $L\geq 0$), then the function $z=B\circ x$ is absolutely continuous, too. This can be easily seen by using the characterization of absolute continuity in
Definition \ref{abs-cont}(iii). Moreover, $z$ is almost everywhere differentiable and the inequality $\|\dot z (\cdot)\|\leq L\|\dot x(\cdot)\|$ holds almost everywhere.
\end{remark}

\section{On the KL property}

The convergence of the trajectories of the dynamical system \eqref{dysy} will be established within the framework of functions satisfying the Kurdyka–Łojasiewicz property. For any \(\eta\in(0,+\infty]\), we denote by \(\Theta _{\eta }\) the class of continuous, concave functions \(\varphi:[0,\eta)\rightarrow [0,+\infty)\) such that $\varphi(0)=0,$ $\varphi$ is continuously differentiable on $(0,\eta)$ with $\varphi'(s)>0$ for all $s\in(0, \eta)$, and $\varphi$ is continuous at $0.$

\begin{definition}\label{KL-property} \rm({\it Kurdyka–Łojasiewicz property}) Let \(f:\R^n\rightarrow \R\) be a differentiable function. We say that $f$ satisfies the {\it Kurdyka–Łojasiewicz (KL) property} at $\ol x\in \R^n$ if there exist \(\eta \in(0,+\infty]\), a neighborhood $U$ of $\ol x,$ and a function \(\varphi\in \Theta_{\eta}\) such that the inequality
$$\varphi'(f(x)-f(\ol x))\Vert{}\n f(x)\Vert{}\geq 1$$
holds for all $x$ belonging to the intersection
$$U\cap \{x\in\R^n: f(\ol x)<f(x)<f(\ol x)+\eta\}.$$
If $f$ satisfies the KL property at every point in \(\mathbb{R}^{n}\), then $f$ is called a {\it KL function}.
\end{definition}
The function $\varphi$ in the definition above is referred to as a desingularizing function \cite{Begout2015}. The origins of this concept trace back to the pioneering work of Łojasiewicz \cite{lojasiewicz1963}, who proved that for any real-analytic function \(f:\R^n\rightarrow\R\) and a critical point $\ol x\in\R^n$ (meaning $\nabla f(\ol x)=0$), there exists an exponent \(\theta\in[1/2,1)\) such that the expression $\vert{}f-f(\ol x)\vert{}^{\theta}\Vert{}\nabla f\Vert{}^{-1}$ remains bounded around $\ol x.$ This corresponds precisely to the setting where \(\varphi(s)=C(1-\theta)^{-1}s^{1-\theta}\). This classical result allows us to interpret the KL property as a re-parametrization of the function values designed to eliminate flatness around critical points.
Kurdyka \cite{kurdyka1998} subsequently extended this property to differentiable functions definable in an o-minimal structure. Further generalizations to the nonsmooth setting can be found in \cite{b-d-l2006, att-b-red-soub2010, b-d-l-s2007, b-d-l-m2010}. One of the most remarkable features of KL functions is their ubiquity in practical applications, see \cite{b-sab-teb, BCL1,ALV-numa,L-mapr,L-jota}. The class of KL functions encompasses semi-algebraic, real sub-analytic, semiconvex, and uniformly convex functions, as well as convex functions satisfying a growth condition. We refer the reader to \cite{b-d-l2006, att-b-red-soub2010, b-d-l-m2010, b-sab-teb, b-d-l-s2007, att-b-sv2013, attouch-bolte2009} and the references therein for comprehensive details and illustrative examples regarding these classes.
A crucial role in our convergence analysis will be played by the following uniformized KL property established in \cite[Lemma 6]{b-sab-teb}. We recall that the distance function to a set \(A\subseteq\R^n\) is defined as $\dist(x,A)=\inf_{y\in A}\Vert{}x-y\Vert{}$ for all \(x\in\R^n\).

\begin{lemma}\label{unif-KL-property} Let \(\Omega\subseteq \R^n\) be a compact set and let \(f:\R^n\rightarrow \R\) be a differentiable function. Assume that $f$ is constant on $\Omega$ and satisfies the KL property at every point of $\Omega.$ Then there exist constants $\e,\eta >0$ and a function \(\varphi\in \Theta_{\eta}\) such that for all $\ol x\in\Omega$ and all $x$ in the intersection
\begin{equation*} \{x\in\R^n: \dist(x,\Omega)<\varepsilon\}\cap \{x\in\R^n: f(\ol x)<f(x)<f(\ol x)+\eta\}
\end{equation*} the following inequality holds:
\begin{equation*}\varphi'(f(x)-f(\ol x))\Vert{}\n f(x)\Vert{}\geq 1.
\end{equation*}
\end{lemma}

Another related notion required for our analysis is the Łojasiewicz exponent, which is defined as follows (see \cite{att-b-red-soub2010,att-b-sv2013,LiP}):
\begin{definition}\label{Lexpo} Let $f:\R^n\To \R$ be a differentiable function. The function $f$ is said to satisfy the Łojasiewicz property if, for every critical point $\ol x\in\crit{f},$ there exist constants $K,\e>0$ and $\t\in(0,1)$ such that
$$\vert{}f(x)-f(\ol x)\vert{}^\t\le K\Vert{}\n f(x)\Vert{}$$
for every $x$ fulfilling $\Vert{}x-\ol x\Vert{}<\e.$ The value $\t$ is called the Łojasiewicz exponent of $f$ at the critical point $\ol x.$ \end{definition}


\begin{thebibliography}{99}


\bibitem{abbas-att-sv} B. Abbas, H. Attouch, B.F. Svaiter, { Newton-like dynamics and forward-backward methods for
structured monotone inclusions in Hilbert spaces}, Journal of Optimization Theory and its Applications 161(2)  (2014) 331-360.

\bibitem{ALV-numa} C.D. Alecsa, S.C. L\'aszl\'o,  A. Viorel, {A gradient-type algorithm with backward inertial steps associated to a nonconvex minimization problem}, Numer Algor 84 (2020) 485–512.

\bibitem{alv-att-bolte-red} F. Alvarez, H. Attouch, J. Bolte, P. Redont, { A second-order gradient-like dissipative dynamical system
with Hessian-driven damping. Application to optimization and mechanics},  Journal de Math\'{e}matiques Pures et Appliqu\'{e}es (9) 81(8) (2002)
747-779.

\bibitem{attouch-bolte2009} H. Attouch, J. Bolte, { On the convergence of the proximal algorithm for nonsmooth functions involving analytic
features}, Mathematical Programming 116(1-2) Series B (2009) 5-16.


\bibitem{att-b-red-soub2010} H. Attouch, J. Bolte, P. Redont, A. Soubeyran, {\it Proximal alternating minimization and projection
methods for nonconvex problems: an approach based on the Kurdyka-\L{}ojasiewicz inequality}, Mathematics of Operations Research
35(2), 438--457, 2010


\bibitem{att-b-sv2013} H. Attouch, J. Bolte, B.F. Svaiter, {\it Convergence of descent methods for semi-algebraic and tame problems:
proximal algorithms, forward-backward splitting, and regularized Gauss-Seidel methods}, Mathematical Programming 137(1-2) Series A, 91-129, 2013

%\bibitem{ABG} H. Attouch, G. Buttazo, G. Michaille, Variational Analysis in Sobolev ans BV Spaces: Applications to PDEs and Optimization, Second Edition, MOS-SIAM Series on Optimization, Phiadelphia, 2014.


\bibitem{att-sv2011} H. Attouch, B.F. Svaiter, {\it A continuous dynamical Newton-like approach to solving monotone inclusions},
SIAM Journal on Control and Optimization 49(2), 574-598, 2011

\bibitem{BauComb} H.H. Bauschke P.L. Combettes, { Convex  Analysis  and  Monotone  Operator  Theory  in Hilbert Spaces}, CMS Books in Mathematics, Springer, New York, 2011

\bibitem{Begout2015} P. B\'egout, J. Bolte, M.A. Jendoubi, {\it On damped second-order gradient systems}, Journal of Differential Equations 259(7) (2015), 3115--3143. DOI: 10.1016/j.jde.2015.04.016

\bibitem{b-sab-teb} J. Bolte, S. Sabach, M. Teboulle, {\it Proximal alternating linearized minimization
for nonconvex and nonsmooth problems}, Mathematical Programming Series A (146)(1--2), 459-494, 2014

\bibitem{b-d-l2006} J. Bolte, A. Daniilidis, A. Lewis, {\it The \L{}ojasiewicz inequality for nonsmooth subanalytic functions with applications
to subgradient dynamical systems}, SIAM Journal on Optimization 17(4), 1205-1223, 2006

\bibitem{b-d-l-s2007} J. Bolte, A. Daniilidis, A. Lewis, M. Shiota, {\it Clarke subgradients of stratifiable functions},
SIAM Journal on Optimization 18(2), 556--572, 2007

\bibitem{b-d-l-m2010} J. Bolte, A. Daniilidis, O. Ley, L. Mazet, {\it Characterizations of \L{}ojasiewicz inequalities:
subgradient flows, talweg, convexity}, Transactions of the American Mathematical Society 362(6), 3319-3363, 2010

\bibitem{BCL} R.I. Bo\c t,  E.R. Csetnek,  S.C. L\' aszl\' o, {\it  Approaching nonsmooth nonconvex minimization through second-order proximal-gradient dynamical systems}, J. Evol. Equ. 18, 1291–1318, 2018

\bibitem{BCL1} R.I. Bo\c t,  E.R. Csetnek,  S.C. L\' aszl\' o, {\it A second-order dynamical approach with variable damping to nonconvex smooth minimization}, Applicable Analysis 99(3), 361–378, 2020



\bibitem{haraux} A. Haraux, {\it Syst\`{e}mes Dynamiques Dissipatifs et Applications},
Recherches en Math\'{e}- matiques Appliqu\'{e}ées 17, Masson, Paris, 1991


\bibitem{kurdyka1998} K. Kurdyka, {\it On gradients of functions definable in o-minimal structures},
Annales de l'institut Fourier (Grenoble) 48(3), 769--783, 1998

\bibitem{L-mapr} S.C. L\' aszl\' o, {\it Convergence rates for an inertial algorithm of gradient type associated to a smooth non-convex minimization}, Math. Program. 190, 285–329, 2021

\bibitem{L-jota} S.C. L\' aszl\' o,  {\it Forward–Backward Algorithm With Different Inertial Terms for Structured Non-Convex Minimization Problems}, J Optim Theory Appl 198, 387–427, 2023
\bibitem{L} S.C. L\'aszl\'o, {\it Solving nonconvex optimization problems via a second-order dynamical system with unbounded damping},
Communications in Optimization Theory 2026, Article 16, 1-19,  2026


\bibitem{LiP}  G. Li, T.K. Pong, {\it Calculus of the exponent of Kurdyka–Łojasiewicz inequality and its applications to linear convergence of first-order methods,} Found. Comput. Math. 18, 1199–1232, 2018

\bibitem{lojasiewicz1963} S. \L{}ojasiewicz, {\it Une propri\'{e}t\'{e} topologique des sous-ensembles analytiques r\'{e}els},
Les \'{E}quations aux D\'{e}riv\'{e}es Partielles, \'{E}ditions du Centre National de la Recherche Scientifique Paris, 87--89, 1963


\bibitem{Nest1} Y. Nesterov, {\it A method of solving a convex programming problem with convergence rate $O(1/k^2)$}, Soviet Math. Dokl.  27, 372-376 1983

\bibitem{rock-wets} R.T. Rockafellar, R.J.-B. Wets, {Variational Analysis}, Fundamental Principles of Mathematical Sciences 317,
Springer-Verlag, Berlin, 1998

\bibitem{sontag} E.D. Sontag, { Mathematical control theory. Deterministic finite-dimensional systems}, Second edition,
Texts in Applied Mathematics 6, Springer-Verlag, New York, 1998.


\bibitem{SBC}W. Su, S. Boyd, E.J. Candes, { A differential equation for modeling Nesterov's
accelerated gradient method: theory and insights}, Journal of Machine Learning Research 17(153) (2016), 1--43.



\end{thebibliography}
\end{document}